\documentclass{IEEEtran}
\usepackage{graphicx}
\usepackage{amsmath,amsthm,amssymb}
\usepackage{subcaption}
\usepackage{gnuplottex}

\usepackage[usenames,dvipsnames,svgnames,table]{xcolor}

\newtheorem{theorem}{Theorem}
\newtheorem{lemma}{Lemma}
\newtheorem{corollary}{Corollary}
\newtheorem{definition}{Definition}

\theoremstyle{definition}
\newtheorem{example}{Example}

\theoremstyle{remark}

\newif\ifexpand
\expandfalse

\def\conv{\mathrm{conv}}

\def\Prob{\mathrm{Prob}}
\def\1{\mathbf{1}}
\newcommand{\vv}[1]{\boldsymbol{#1}}
\newcommand{\vvi} [2]{ \boldsymbol{#1}^{(#2)} }
\newcommand{\vi} [2]{ {#1}^{(#2)} }

\title{Max-Weight Revisited: Sequences of Non-Convex Optimisations Solving Convex Optimisations}
\author{V\'ictor Valls, Douglas J. Leith\\Trinity College Dublin\thanks{This work
  was supported by Science Foundation Ireland
  under Grant No. 11/PI/1177.}}

\begin{document}
\maketitle

\begin{abstract}
We investigate the connections between max-weight approaches and dual subgradient methods for convex optimisation.  We find that strong connections exist and we establish a clean, unifying theoretical framework that includes both max-weight and dual subgradient approaches as special cases. Our analysis uses only elementary methods, and is not asymptotic in nature. It also allows us to establish an explicit and direct connection between discrete queue occupancies and Lagrange multipliers. 
\end{abstract}


\section{Introduction}\label{sec:introduction}
\IEEEPARstart{I}{n} queueing networks, max-weight (also referred to as backpressure) approaches have been the subject of much interest for solving utility optimisation problems in a distributed manner.

In brief, consider a queueing network where the queue occupancy of the $i$'th queue at time $k$ is denoted by $\vi{Q}{i}_{k}\in\mathbb{N}$, $i=1,2,\dots,n$, and we gather these together into vector $\vv{Q}_k\in\mathbb{N}^n$.  Time is slotted and at each time step  $k=1,2,\dots$ we select action $\vv{x}_k\in D\subset\mathbb{N}^n$, \emph{e.g.,} selecting $i$'th element $\vi{x}{i}_k=1$ corresponds to transmitting one packet from queue $i$ and $\vi{x}{i}_k=0$ to doing nothing.   The connectivity between queues is captured via matrix $\vv{A}\in\{-1,0,1\}^{n\times n}$, whose $i$'th row has a $-1$ at the $i$'th entry, $1$ at entries corresponding to queues from which packets are sent to queue $i$, and $0$ entries elsewhere. The queue occupancy then updates according to $\vv{Q}_{k+1}=[\vv{Q}_k+\vv{A}\vv{x}_k+\vv{b}_k]^+$, $i=1,2,\dots,n$, where the $i$'th element of vector $\vv{b}_k\in\mathbb{N}^n$ denotes the number of external packet arrivals to queue $i$ at time $k$.  The objective is to stabilise all of the queues while maximising utility ${U(\vv{z}_k)}$ where $U:\mathbb{R}^n\rightarrow\mathbb{R}$ is concave and continuously differentiable and $\vv{z}_k$ is a running average of $\vv{x}_j$, $j=1,\dots,k$.  The \emph{greedy primal-dual} variant of max-weight scheduling \cite{stolyargreedy}, for example, selects action $\vv{x}_k \in \arg\max_{\vv{x}\in D} \partial U(\vv{z}_k)^T \vv{x}- \beta \vv{Q}^T_k\vv{A}\vv{x}$ with $\vv{z}_{k+1}=(1-\beta)\vv{z}_k+\beta\vv{x}_k$, $0<\beta <1$ a design parameter. 

Appealing features of this max-weight scheduling approach include the lack of a requirement for \emph{a priori} knowledge of packet arrival process $\{\vv{b}_k\}$, and the fact that the discrete action set matches the actual decision variables (namely, do we transmit a packet or not).  Importantly, although cost function $-U(\cdot)$ is required to be convex, at each time step the max-weight optimisation is non-convex owing to the non-convexity of action set $D$.  Further, convergence is typically proved using Foster-Lyapunov or by sophisticated fluid-limit arguments, which allow sequence $\{\vv{b}_k\}$ to be accommodated but are distinct from the usual approaches employed in convex optimisation.   Hence, the body of work on max-weight approaches remains separate from the mainstream literature on convex optimisation.    On the other hand, queueing and Lagrange multiplier subgradient updates are clearly similar, at least superficially, although the exact nature of the relationship between queues and multipliers remains unclear. 

Taking these observations as our starting point, in this paper we investigate the connections between max-weight approaches and dual subgradient methods for convex optimisation.  We find that strong connections do indeed exist and we establish a clean, unifying theoretical framework that includes both max-weight and dual subgradient approaches as special cases.
In summary, the main contributions of the paper include the following.   


1) \emph{Generalising max-weight}.  Our analysis places max-weight firmly within the field of convex optimisation, extending it from the specific constraints induced by queueing networks to general convex nonlinear contraints with bounded curvature.   We show that any non-convex update with suitable descent properties can be employed, and the wealth of convex descent methods can be leveraged to derive non-convex approaches.  Descent methods studied here include non-convex variants of the classical Frank-Wolfe update and of the primal Lagrangian update.


2) \emph{Generalising dual subgradient methods}.  We show that convexity can be relaxed in classical dual subgradient methods, allowing use of a finite action set.   In the special case of optimisation problems with linear constraints, we rigorously establish a close connection (essentially an equivalence) between Lagrange multiplier subgradient updates and discrete queues, so putting existing intuition on a sound footing.   

3) \emph{Unifying theoretical framework}.  In generalising max-weight and dual subgradient methods our analysis clarifies the fundamental properties required.   In particular, bounded curvature of the objective and constraint functions plays a prominent role in our analysis, as does boundedness of the action set.  Of interest in its own right, we note that our analysis requires only elementary methods and so an additional contribution is the accessible nature of the methods of proof employed. In particular, it turns out that deterministic analysis of sample paths is sufficient to handle stochasticity. The methods of proof themselves are new in the context of max-weight approaches, and are neither Foster-Lyapunov nor fluid-limit based.    



\subsection{Related Work}

Max-weight scheduling was introduced by Tassiulas and Ephremides in their seminal paper \cite{tassiulasstability}.  They consider a network of queues with slotted time, an integer number of packet arrivals in each slot and a finite set of admissible scheduling patterns, referred to as \emph{actions}, in each slot.   Using a Forster-Lyapunov approach they present a scheduling policy that stabilises the queues provided the external traffic arrivals are strictly feasible. Namely, the scheduling policy consists of selecting the action at each slot that maximises the queue-length-weighted sum of rates, $\vv{x}_k \in \arg \max_{\vv{x} \in D} -  \vv{Q}_k^T \vv{A}\vv{x} $.  

Independently, \cite{stolyargreedy,neelypower,eryilmazfair} proposed extensions to the max-weight approach to accommodate concave utility functions.   
In \cite{stolyargreedy} the \emph{greedy primal-dual} algorithm is introduced, as already described above, for network linear constraints and utility function $U(\cdot)$ which is continuously differentiable and concave. The previous work is extended in \cite{stolyar2006greedy} to consider general nonlinear constraints.
In \cite{eryilmazfair} the utility fair allocation of throughput in a cellular downlink is considered.   The utility function is of the form $U(\vv{z})=\sum_{i=1}^n U_i(\vi{z}{i})$, $U_i(z)=\beta_i ({z^{(1-\frac{1}{m})}})/({1-\frac{1}{m}})$, with $m$, $\beta_i$ design parameters.   Queue departures are scheduled according to  $\vv{x}_k \in \arg \max_{\vv{x} \in \conv(D)} -  \vv{Q}_k^T \vv{A}\vv{x}$ and queue arrivals are scheduled by a congestion controller such that $E[\vi{b}{i}_k|\vv{Q}_k] = \min\{\partial U_i(\vi{Q}{i}_k),M\}$ and  $E[(\vi{b}{i}_k)^2|\vv{Q}_k]\le A$ where $A$, $M$ are positive constants. 
The work in \cite{neelypower} considers power allocation in a multibeam downlink satellite communication link with the aim of maximising throughput while ensuring queue stability.  This is extended in a sequence of papers \cite{neelydynamic,neelyenergy, neelyfairness} and a book \cite{neelybook} to develop the \textit{drift plus penalty} approach. In this approach the basic strategy for scheduling queue departures is according to $\vv{x}_k \in \arg \max_{\vv{x} \in D} -  \vv{Q}_k^T \vv{A}\vv{x}$ and utility functions are incorporated in a variety of ways.  For example, for concave non-decreasing {continuous} utility functions $U$ of the form $U(\vv{z})=\sum_{i=1}^nU_i(\vi{z}{i})$ one formulation is for a congestion controller to schedule arrivals into an ingress queue such that $\vi{b}{i}_k\in \arg\max_{0\le b \le R}VU_i(b)-b\vi{Q}{i}_k$ where $V$, $R$ are sufficiently large design parameters and $b\in\mathbb{R}$ \cite{tassiulas_neely_book}.  Another example is for cost functions of the form $E[P_k(\vv{x}_k)]$ where $P_k(\cdot)$ is bounded, i.i.d. and known at each time step,  in which case actions at each time step are selected to minimise $\vv{x}_k \in \arg \min_{\vv{x} \in D} VP_k(\vv{x}_k)+\vv{Q}_k^T \vv{A}\vv{x}$ where $V$ is a design parameter \cite{neelybook}.

With regard to the existence of a connection between the discrete-valued queue occupancy in a queueing network and continuous-valued Lagrange multipliers, this has been noted by several authors, see for example \cite{neely05,lin06}, and so might be considered something of a ``folk theorem'' but we are aware of few rigorous results.  A notable exception is \cite{neely11}, which establishes that a discrete queue update tends on average to drift towards the optimal multiplier value. Also, the \emph{greedy primal-dual} algorithm presented in \cite{stolyargreedy} shows that asymptotically as design parameter $\beta \rightarrow 0$ and $t \rightarrow \infty$ the scaled queue occupancy converges to the set of dual optima.

Selection of a sequence of actions in a discrete-like manner is also considered in the convex optimisation literature.  The \emph{nonlinear Gauss-Seidel} algorithm, also known as \emph{block coordinate descent} \cite{bertsekasnonlinear,bertsekasparallel} minimises a convex function over a convex set by updating one co-ordinate at a time.  The convex function is required to be continuously differentiable and strictly convex and, unlike in the max-weight algorithms discussed above, the action set is convex.  The classical Frank-Wolfe algorithm \cite{frankwolfe} also minimises a convex continuously differentiable function over a polytope by selecting from a discrete set of descent directions, although a continuous-valued line search is used to determine the final update.   We also note the work on online convex optimisation \cite{zinkevichonline,flaxmanonline}, where the task is to choose a sequence of actions so to minimise an unknown sequence of convex functions with low regret.  


\subsection{Notation}
Vectors and matrices are indicated in bold type.  Since we often use subscripts to indicate elements in a sequence, to avoid confusion we usually use a superscript $\vi{x}{j}$ to denote the $j$'th element of a vector $\vv{x}$.  The $j$'th element of operator $[\vv{x}]^{[0,\bar{\lambda}]}$ equals $\vi{x}{j}$ (the $j$'th element of $\vv{x}$) when $\vi{x}{j}\in{[0,\bar{\lambda}]}$ and otherwise equals $0$ when $\vi{x}{j}<0$ and $\bar{\lambda}$ when $\vi{x}{j} > \bar{\lambda}$.  Note that we allow $\bar{\lambda}=+\infty$, and following standard notation in this case usually write $[x]^+$ instead of $[x]^{[0,\infty)}$.  The subgradient of a convex function $f$ at point $\vv{x}$ is denoted $\partial f(\vv{x})$. {For two vectors $\vv{x},\vv{y} \in \mathbb R^m$ we use element-wise comparisons $\vv{x} \succeq \vv{y}$ and $\vv{y} \succ \vv{x}$ to denote when $\vi{y}{j} \ge \vi{x}{j}$, $\vi{y}{j} > \vi{x}{j}$ respectively for all $j=1,\dots,m$.}
\section{Preliminaries}

We recall the following convexity properties.

\begin{lemma} [Lipschitz Continuity]\label{lem:lipschtitz}
Let $h:M\rightarrow \mathbb{R}$ be a convex function and let $C$ be a closed and bounded set contained in the relative interior of the domain $M\subseteq \mathbb{R}^n$. Then $h(\cdot)$ is Lipschitz continuous on $C$ \emph{i.e.,} there exists constant $\nu_h$  such that $|h(\vv{x}) - h(\vv{y})| \le \nu_h \|\vv{x} -\vv{y}\|_2$ $\forall \vv{x}, \vv{y} \in C$. 
\end{lemma}
\begin{IEEEproof}
 See, for example, \cite{lipschitz}.
\end{IEEEproof}

\begin{lemma}[Bounded Distance]\label{lem:one}
Let  $D:=\{\vv{x}_{1},\dots,\vv{x}_{|D|}\}$ be a finite set of points from $\mathbb{R}^n$.  Then there exists constant $\bar{x}_D$ such that $\|\vv{z}-\vv{y}\|_2 \le \bar{x}_D$ for any two points $\vv{z},\vv{y}\in C:=\conv(D)$, where $\conv(D)$ denotes the convex hull of $D$.
\end{lemma}

\begin{IEEEproof}
Since $\vv{z},\vv{y}\in C$ these can be written as the convex combination of points in $D$, \emph{i.e.,} $\vv{z}=\sum_{j=1}^{|D|}\vi{a}{j} \vv{x}_j$, $\vv{y}=\sum_{j=1}^{|D|}\vi{b}{j} \vv{x}_j$ with $\| \vv{a} \|_1 = 1 = \| \vv{b}\|_1$. Hence $\|\vv{z}-\vv{y}\|_2 = \|\sum_{j=1}^{|D|}(\vi{a}{j}-\vi{b}{j})\vv{x}_j\|_2 \le \sum_{j=1}^{|D|}\| \vi{a}{j}-\vi{b}{j}\|_2 \|\vv{x}_j\|_2 \le \bar{x}_D:=2 \max_{\vv{x}\in D}\|\vv{x}\|_2$.
 \end{IEEEproof}
 
We also introduce the following definition:
\begin{definition}[Bounded Curvature]
Let $h:M\rightarrow\mathbb{R}$ be  a convex function defined on domain $M\subseteq \mathbb{R}^n$.   We say the $h(\cdot)$ has bounded curvature on set $C\subset M$ if for any points $\vv{z}, \vv{z}+\vv{\delta} \in C$
\begin{align}\label{eq:bound}
h(\vv{z}+\vv{\delta}) - h(\vv{z}) \le    \partial h(\vv{z})^T \vv{\delta} + \mu_h \| \vv{\delta} \|_2^2
\end{align}
where {$\mu_h \ge 0$} is a constant that does not depend on $\vv{z}$ or $\vv{\delta}$.  
\end{definition}
Bounded curvature will prove important in our analysis.  The following lemma shows that a necessary and sufficient condition for bounded curvature is that the subgradients of $h(\cdot)$ are Lipschitz continuous on set $C$. 
\begin{lemma}[Bounded Curvature] \label{lem:curv}
Let $h:M\rightarrow\mathbb{R}$, $M\subseteq \mathbb{R}^n$ be  a convex function.  Then $h(\cdot)$ has bounded curvature on $C$ if and only if for all $\vv{z}, \vv{z}+\vv{\delta} \in C$ there exists a member $\partial h(\vv{z})$ (respectively, $\partial h(\vv{z}+\vv{\delta})$) of the set of subdifferentials at point $\vv{z}$ (respectively, $\vv{z}+\vv{\delta}$) such that $\left(\partial h(\vv{z}+\vv{\delta})-\partial h(\vv{z})\right)^T\vv{\delta}  \le \mu_h \| \vv{\delta} \|_2^2$ where $\mu_h$ does not depend on $\vv{z}$ or $\vv{\delta}$.
\end{lemma}
\begin{IEEEproof}
$\Rightarrow$ Suppose $h(\cdot)$ has bounded curvature on $C$.  From (\ref{eq:bound}) it follows that $h(\vv{z}+\vv{\delta}) - h(\vv{z}) \le    \partial h(\vv{z})^T \vv{\delta} + \mu_h \| \vv{\delta} \|_2^2$ and $h(\vv{z}) - h(\vv{z}+\vv{\delta}) \le    -\partial h(\vv{z}+\vv{\delta})^T \vv{\delta} + \mu_h \| \vv{\delta} \|_2^2$.  Adding left-hand and right-hand sides of these inequalities yields $0 \le \left(\partial h(\vv{z}) -\partial h(\vv{z}+\vv{\delta})\right)^T \vv{\delta} + 2\mu_h \| \vv{\delta} \|_2^2$ \emph{i.e.,} $\left(\partial h(\vv{z}+\vv{\delta})-\partial h(\vv{z})\right)^T\vv{\delta}  \le \mu_h \| \vv{\delta} \|_2^2$.  

$\Leftarrow$ Suppose $\left(\partial h(\vv{z}+\vv{\delta}) - \partial h(\vv{z}) \right)^T\vv{\delta} \le \mu_h \|\vv{\delta} \|_2$ for all $\vv{z}, \vv{z}+\vv{\delta} \in M$. It follows that $\partial h(\vv{z}+\vv{\delta})^T\vv{\delta} \le \partial h(\vv{z})^T \vv{\delta} + \mu_h \|\vv{\delta} \|_2^2$.  By the definition of the subgradient we have that $h(\vv{z}+\vv{\delta})-h(\vv{z}) \le \partial h(\vv{z}+\vv{\delta})^T \vv{\delta}$, and so we obtain that $h(\vv{z}+\vv{\delta}) - h(\vv{z}) \le    \partial h(\vv{z})^T \vv{\delta} + \mu_h \| \vv{\delta} \|_2^2.$
\end{IEEEproof}

\section{Non-Convex Descent}\label{sec:nonconvexdescent}
 
We begin by considering minimisation of convex function $F:\mathbb{R}^n\rightarrow\mathbb{R}$ on convex set $C:=\conv(D)$, the convex hull of set $D:=\{\vv{x}_{1},\dots,\vv{x}_{|D|}\}$ consisting of a finite collection of points from $\mathbb{R}^n$ (so $C$ is a polytope).   Our interest is in selecting a sequence of points $\{\vv{x}_k\}$, $k=1,2,\dots$ from set $D$ such that the running average $\vv{z}_{k+1}=(1-\beta)\vv{z}_k+\beta \vv{x}_k$ minimises $F(\cdot)$ for $k$  sufficiently large and $\beta$ sufficiently small.    Note that set $D$ is non-convex since it consists of a finite number of points, and by analogy with max-weight terminology we will refer to it as the \emph{action set}.


Since $C$ is the convex hull of action set $D$, any point $\vv{z}^*\in C$ minimising $F(\cdot)$ can be written as convex combinations of points in $D$ \emph{i.e.,} $\vv{z}^*=\sum_{j=1}^{|D|}\vi{\theta^*}{j} \vv{x}_j$, $\vi{\theta^*}{j}\in[0,1]$, $\| \vv{\theta} \|_1=1$.   Hence, we can always construct sequence $\{\vv{x}_k\}$ by selecting points from set $D$ in proportion to the $\vi{\theta^*}{j}$, $j=1,\dots,|D|$.  That is, by \emph{a posteriori} time-sharing (\emph{a posteriori} in the sense that we need to find minimum $\vv{z}^*$ before we can construct sequence  $\{\vv{x}_k\}$). Of more interest, however, it turns out that when function $F(\cdot)$ has bounded curvature then sequences $\{\vv{x}_k\}$ can be found without requiring knowledge of $\vv{z}^*$.   

\subsection{Non-Convex Direct Descent}
The following theorem formalises the above commentary, also generalising it to sequences of convex functions $\{F_k\}$ rather than just a single function as this will prove useful later.

 \begin{theorem}[Greedy Non-Convex Convergence]\label{th:descent}
Let $\{F_k\}$ be a sequence of convex functions with uniformly bounded curvature $\mu_F$ on set $C:=\conv(D)$, action set $D$ a finite set of points from $\mathbb{R}^n$.   Let $\{\vv{z}_k\}$ be a sequence of vectors satisfying $\vv{z}_{k+1}=(1-\beta)\vv{z}_k+\beta \vv{x}_k$ with $\vv{z}_1\in C$ and
\begin{align}
\vv{x}_k\in\arg\min_{\vv{x}\in D} F_k((1-\beta)\vv{z}_k+\beta \vv{x}),\ k=1,2,\dots \label{eq:th1}
\end{align}
 Suppose parameter $\beta$ is sufficiently small that 
 \begin{align}
{0 <} \beta\le(1-\gamma)\min\{\epsilon / (\mu_F\bar{x}^2_D),1\} \label{eq:betath1}
 \end{align}
 with $\epsilon>0$, $\gamma \in (0,1)$, $\bar{x}_D:=2 \max_{\vv{x}\in D}\|\vv{x}\|_2$ and that functions $F_k$ change sufficiently slowly that 
 \begin{align*}
 |F_{k+1}(\vv{z})-F_k(\vv{z})|\le \gamma_1 \gamma \beta \epsilon,\  \forall \vv{z}\in C
 \end{align*}
 with $ \gamma_1 \in (0,{1}/{2})$.   Then for every $\epsilon >0$ and $k$ sufficiently large we have that
\begin{align*}
0\le  F_k(\vv{z}_{k+1})-F_k(\vv{y}_k^*) \le 2\epsilon
 \end{align*}
 where $\vv{y}_k^* \in \arg\min_{\vv{z}\in C} F_k(\vv{z})$. 
\end{theorem}
\begin{IEEEproof} 
See Appendix.
\end{IEEEproof}

Observe that in Theorem \ref{th:descent} we select $\vv{x}_k$ by solving non-convex optimisation (\ref{eq:th1}) at each time step.   This optimisation is one step ahead, or greedy, in nature and does not look ahead to future values of the sequence or require knowledge of optima $\vv{y}^*_k$.   Of course, such an approach is mainly of interest when non-convex optimisation (\ref{eq:th1}) can be efficiently solved, \emph{e.g.,} when action set $D$ is small or the optimisation separable.   

\begin{figure}
\centering
\begin{subfigure}[b]{0.4\textwidth}
\centering
\includegraphics[width=0.85\textwidth]{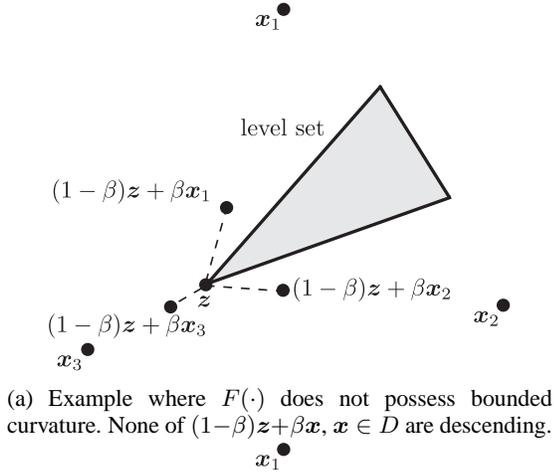}
\caption{Example where $F(\cdot)$ does not possess bounded curvature.  None of $(1-\beta)\vv{z}+\beta \vv{x}$, $\vv{x}\in D$ are descending. }\label{fig:ex}
\end{subfigure}
\qquad
\begin{subfigure}[b]{0.4\textwidth}
\centering
\includegraphics[width=0.85\textwidth]{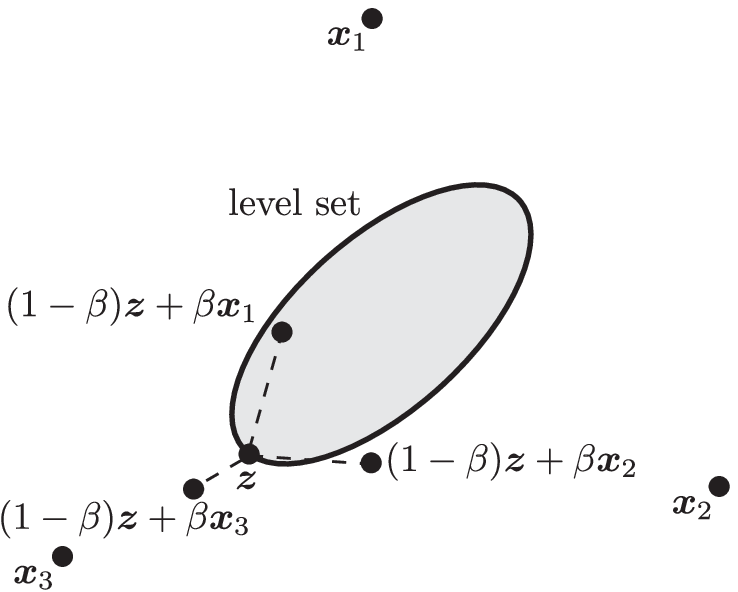}
\caption{Example where $F(\cdot)$ has bounded curvature.  For $\beta$ sufficiently small, for at least one $(1-\beta)\vv{z}+\beta \vv{x}$, $\vv{x}\in D$ descent is possible.}\label{fig:ex2}
\end{subfigure}
\caption{Illustrating how bounded curvature allows monotonic descent.  Set $D$ consists of the marked points $\vv{x}_{1}$, $\vv{x}_{2}$, $\vv{x}_{3}$.  Level set $\{F(\vv{y}) \le F(\vv{z}):\vv{y}\in C\}$ is indicated by the shaded areas.  The possible choices of $(1-\beta)\vv{z}+\beta \vv{x}$, $\vv{x}\in D$ are indicated.}
\end{figure}

Observe also that Theorem \ref{th:descent} relies upon the bounded curvature of the sequence of functions $F_k(\cdot)$.     A smoothness assumption of this sort seems essential, since when it does not hold it is easy to construct examples where Theorem \ref{th:descent} does not hold.   Such an example is illustrated schematically in Figure \ref{fig:ex}.    The shaded region in Figure \ref{fig:ex} indicates the level set $\{F(\vv{y}) \le F(\vv{z}):\vv{y}\in C\}$.  The level set is convex, but the boundary is non-smooth and contains ``kinks''.  We can select points from the set $\{(1-\beta) \vv{z} + \beta \vv{x}: \vv{x}\in D=\{\vv{x}_{1},\vv{x}_{2},\vv{x}_{3}\}\}$.   This set of points is indicated in Figure \ref{fig:ex} and it can be seen that  every point  lies outside the level set.  Hence, we must have $F((1-\beta) \vv{z} + \beta \vv{x}) > F(\vv{z})$, and upon iterating we will end up with a diverging sequence.  Note that in this example changing the step size $\beta$ does not resolve the issue.  Bounded curvature ensures that the boundary of the level sets is smooth, and this ensures that for sufficiently small $\beta$ there exists a convex combination of $\vv{z}$ with a point $\vv{x} \in D$ such that $F((1-\beta) \vv{z} + \beta \vv{x}) < F(\vv{z})$ and so the solution to optimisation (\ref{eq:th1}) improves our objective, see Figure \ref{fig:ex2}.

Theorem \ref{th:descent} is stated in a fairly general manner since this will be needed for our later analysis.   An immediate  corollary to  Theorem \ref{th:descent} is the following convergence result for unconstrained optimisation.  

\begin{figure}
\centering
\includegraphics[width=0.9\columnwidth]{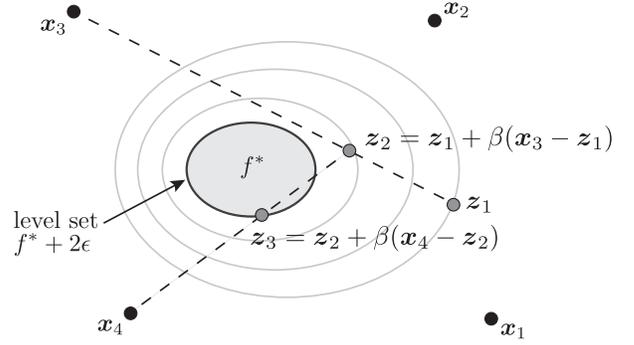}
\caption{Illustrating unconstrained convergence in $\mathbb R^2$. The sequence of non-convex optimisations converges with $k=2$. The function average decreases monotonically and then remains in level set $f(\vv{z}_k) \le f^* + 2 \epsilon$  for $k \ge 3$. }
\label{fig:converence}
\end{figure}

\begin{corollary}[Unconstrained Optimisation]\label{th:unconstrained}
Consider the following sequence of non-convex optimisations $\{P^u_k\}$:
\begin{align*}
&\vv{x}_k \in \arg \underset{\vv{x}\in D}{\min }  f\left((1-\beta)\vv{z}_{k}+\beta \vv{x}\right) \\
&\vv{z}_{k+1}= (1-\beta)\vv{z}_{k}+\beta \vv{x}_k
\end{align*}
with $\vv{z}_1\in C:=\conv{(D)}$, action set $D\subset\mathbb{R}^n$ finite.    Then  $0\le f(\vv{z}_k) - f^* \le 2\epsilon$ for all $k$ sufficiently large, where $f^*=\min_{\vv{z}\in C} f(\vv{z})$, provided $f(\cdot)$ has bounded curvature with curvature constant $\mu_f$ and $0<\beta \le (1-\gamma)\min \{\epsilon / (\mu_f\bar{x}_D^2),1 \}$ with $\gamma \in (0,1)$, $\epsilon>0$, $\bar{x}_D:=2 \max_{\vv{x}\in D}\|\vv{x}\|_2$.
\end{corollary}
Figure \ref{fig:converence} illustrates Corollary \ref{th:unconstrained} schematically in $\mathbb R^2$. The sequence of non-convex optimisations descends in two iterations $f(\vv{z}_1) > f(\vv{z}_2) > f(\vv{z}_3)$ (using points $\vv{x}_3$ and $\vv{x}_4$ respectively) and $f(\vv{z}_k) - f^* \le 2 \epsilon$ for $k > 3$ (not shown in Figure \ref{fig:converence}).  

Note that the curvature constant $\mu_f$ of function $f$ need not be known, an upper bound being sufficient to select $\beta$. Next we present two brief examples that are affected differently by constant $\mu_f$. 

\begin{example}[Linear Objective] \label{ex:linear}Suppose $f (\vv{z}): = \vv{a}^T \vv{z}$ where $\vv{a} \in \mathbb R^n$. The objective function is linear and so has curvature constant $\mu_f = 0$.  It can be seen from (\ref{eq:betath1}) that we can choose $\beta$ independently of parameter $\epsilon$.  Further, for any $\beta \in (0,1)$ we have that $f(\vv{z}_{k+1}) < f(\vv{z}_k)$ for all $k$ and so $f(\vv{z}_k)\rightarrow f^*$.  
\label{ex:1}
\end{example}

\begin{example}[Quadratic Objective] \label{ex:quadratic} Suppose $f(\vv{z}) : = \frac{1}{2}\vv{z}^T \vv{A} \vv{z}$ where $\vv{A} \in \mathbb R^{n\times n}$ is positive definit.  Then $\mu_f = \lambda_{\max} (\vv{A}) > 0$ and in contrast to Example \ref{ex:1} the bound (\ref{eq:betath1}) on parameter $\beta$ now depends on $\epsilon$ and convergence is into the ball $f(\vv{z}_k) - f^* \le 2 \epsilon$ for $k$ sufficiently large.
\label{ex:2}
\end{example}


\subsection{Non-Convex Frank-Wolfe-like Descent}
It is important to note that other convergent non-convex updates are also possible.  For example:
\begin{theorem}[Greedy Non-Convex FW Convergence]\label{th:descent1b}
Consider the setup in Theorem \ref{th:descent}, but with modified update
\begin{align}
\vv{x}_k\in\arg\min_{\vv{x}\in D} \partial F_k(\vv{z}_k)^T\vv{x},\ k=1,2,\dots \label{eq:th1b}
\end{align}
Then for every $\epsilon > 0$ and $k$ sufficiently large we have that
\begin{align*}
0\le  F_k(\vv{z}_{k+1})-F_k(\vv{y}_k^*) \le 2\epsilon
 \end{align*}
 where $\vv{y}_k^* \in \arg\min_{\vv{z}\in C} F_k(\vv{z})$. 
\end{theorem}
\begin{IEEEproof}
See Appendix.
\end{IEEEproof}
The intuition behind the update in Theorem \ref{th:descent1b} is that at each step we locally approximate $F_k(\vv{z}_{k+1})$ by linear function $F_k(\vv{z}_{k})+\partial F_k(\vv{z}_k)^T(\vv{z}_{k+1}-\vv{z}_k)$ and then minimise this linear function.  Since $F_k(\cdot)$ is convex, this linear function is in fact the supporting hyperplane to $F_k(\cdot)$ at point $\vv{z}_k$, and so can be expected to allow us to find a descent direction.   Similar intuition also underlies classical Frank-Wolfe algorithms for convex optimisation \cite{frankwolfe} on a polytope, and Theorem \ref{th:descent1b} extends this class of algorithms to make use of non-convex update (\ref{eq:th1b}) and a fixed step size (rather than the classical approach of selecting the step size by line search).

Note that when the function is linear $F_k(\vv{z})=\vv{c}_k^T\vv{z}$, $\vv{c}_k\in\mathbb{R}^n$, {then $\arg\min_{\vv{x}\in D} F_k((1-\beta) \vv{z} + \beta \vv{x}) = \arg\min_{\vv{x}\in D} \vv{c}_k^T\vv{x}$} and $\arg\min_{\vv{x}\in D} \partial F_k(\vv{z}_k)^T \vv{x}=\arg\min_{\vv{x}\in D}\vv{c}_k^T\vv{x}$.  That is, updates (\ref{eq:th1}) and (\ref{eq:th1b}) are identical.

Note also that 
\begin{align}
\arg\min_{\vv{x}\in D} \partial F_k(\vv{z}_k)^T\vv{x}
\subseteq\arg\min_{\vv{z}\in C} \partial F_k(\vv{z}_k)^T\vv{z}\label{eq:linear}.
\end{align}   
This is because the RHS of (\ref{eq:linear}) is a linear program (the objective is linear and set $C$ is a polytope, so defined by linear constraints) and so the optimum set is either (i) an extreme point of $C$ and so a member of set $D$, or (ii) a face of polytope $C$ with the extreme points of the face belonging to set $D$.   Hence, while update (\ref{eq:th1b}) is non-convex it can nevertheless be solved in polynomial time.

\section{Sequences of Non-Convex Optimisations \& Constrained Convex Optimisation} \label{sec:convc}

We now extend consideration to the constrained convex optimisation $P$:  
\begin{align*}
& \underset{\vv{z}\in C}{\text{minimise  }} \qquad  f(\vv{z}) \\
& \underset{}{\text{subject to}} \qquad  \vv{g}(\vv{z}) \preceq 0
 \end{align*}
 where $\vv{g}(\vv{z}) := [\vi{g}{1},\dots,\vi{g}{m}]^T$ and $f,{\vi{g}{j}}:\mathbb{R}^n\rightarrow \mathbb{R}$, $j=1,\dots,m$ are convex functions with bounded curvature with, respectively, curvature constants $\mu_f$ and $\mu_{{\vi{g}{j}}}$. As before, action set $D$ consists of a finite set of points in $\mathbb R^n$ and $C:=\conv{(D)}$.    Let $C_0:=\{\vv{z}\in C \mid \vv{g} (\vv{z})\preceq \vv{0} \}$ denote the set of feasible points, which we will assume has non-empty relative interior (\emph{i.e.,} a Slater point exists).   Let $C^*:= \arg\min_{\vv{z}\in C_0} f(\vv{z})\subseteq C_0 $ be the set of optima and $f^*:=f(\vv{z}^*)$, $\vv{z}^*\in C^*$. 
 
In the next sections we introduce a generalised dual subgradient approach for finding approximate solutions to optimisation $P$ which, as we will see, includes the classical convex dual subgradient method as a special case. 

\subsection{Lagrangian Penalty}
As in classical convex optimisation we define Lagrangian $L(\vv{z}, \vv{\lambda}) := f(\vv{z}) + \vv{\lambda}^T \vv{g}(\vv{z})$ where $\vv{\lambda} = [\vi{\lambda}{1}, \dots,\vi{\lambda}{m}]^T$ with $\vi{\lambda}{j} \ge 0$, $j=1,\dots,m$.   Since set $C_0$ has non-empty relative interior, the Slater condition is satisfied and strong duality holds. That is, there is zero duality gap and so the solution of the dual problem $P^D$: 
\begin{align*}
\underset{\vv{\lambda} \succeq \vv{0}}{\text{maximise}} \quad q(\vv{\lambda}) := \min_{\vv{z} \in C} L(\vv{z}, \vv{\lambda})
\end{align*}and primal problem $P$ coincide. Therefore, we have that
\begin{align*}
\min_{\vv{z}\in C}\max_{\vv{\lambda}\succeq 0} L(\vv{z},\vv{\lambda}) =\max_{\vv{\lambda}\succeq 0}\min_{\vv{z}\in C} L(\vv{z},\vv{\lambda})=q(\vv{\lambda}^*) = f^*
\end{align*}where $\vv{\lambda}^*: = \arg \max_{\vv{\lambda} \succeq \vv{0}} q(\vv{\lambda})$.

\subsubsection{Lagrangian Bounded Curvature}
As already noted, bounded curvature plays a key role in ensuring convergence to an optimum when selecting from a discrete set of actions.   For any two points $\vv{z},\vv{z}+\vv{\delta} \in C$ we have that
\begin{align*}
L(\vv{z}+\vv{\delta},\vv{\lambda}) \le L(\vv{z},\vv{\lambda}) + \partial_{\vv{z}} L (\vv{z}, \vv{\lambda})^T \vv{\delta} + \mu_L \|\vv{\delta} \|_2 ^2,
\end{align*}
where $
\mu_L = \mu_f + \vv{\lambda}^T \vv{\mu}_{\vv{g}}$ with $\vv{\mu}_{\vv{g}} := [\mu_{\vi{g}{1}}, \dots, \mu_{\vi{g}{m}}]^T$.  It can be seen that the curvature constant $\mu_L$ of the Lagrangian depends on the multiplier $\vv{\lambda}$.  Since set $\vv{\lambda}\succeq \vv{0}$ is unbounded, it follows that the Lagrangian does not have bounded curvature on this set unless $\vv{\mu}_{\vv{g}}=\vv{0}$ (corresponding to the special case where the constraints are linear). Fortunately, by constraining $\vi{\lambda}{j} \le \bar \lambda, \ j=1,\dots,m$ for some $\bar \lambda \ge 0$ resolves the issue, \emph{i.e.,} now $L(\cdot,\vv{\lambda})$ has uniform bounded curvature with constant
\begin{align*}
\bar \mu_L = \mu_f + \bar \lambda \vv{1}^T \vv{\mu}_{\vv{g}}.
\end{align*}
For bounded curvature we only require constant $\bar\lambda$ to be finite, but as we will see later in Lemmas \ref{th:subgradient} and \ref{th:slackness} in general  it should be chosen with some care. 

\subsection{Non-Convex Dual Subgradient Update}
In this section we present a {primal-dual}-like approach in which we use discrete actions to obtain {approximate} solutions to problem $P$. {\color{black}{In particular, we construct a sequence $\{ \vv{z}_k \}$ of points in $C$ such that $f( \frac{1}{k} \sum_{i=1}^k\vv{z}_{i+1})$ is arbitrarily close to $f^*$ for $k$ sufficiently large.  }}


We start by introducing two lemmas, which will play a prominent role in later proofs.


\begin{lemma}[Minimising Sequence of Lagrangians] \label{th:lup} {\color{black}{ Let $\{{\vv{\lambda}}_k \}$ be a sequence of vectors in $\mathbb R^m$ such that ${\vv{\lambda}}_k \preceq \bar \lambda \vv{1}$, $\bar \lambda > 0$ and $\| {\vv{\lambda}}_{k+1} - {\vv{\lambda}}_k \|_2 \le \gamma_1 \gamma \beta \epsilon / ({m} {\bar g})$ with $\gamma \in (0,1)$, $\gamma_1 \in (0, 1/2)$, $\beta$, $\epsilon >0$, ${{\bar g}}:= \max_{\vv{z}\in C} \|\vv{g}(\vv{z})\|_\infty$.
Consider optimisation problem $P$ and updates 
\begin{align}
\vv{x}_{k} & \in \arg \min_{\vv{x} \in D} L((1-\beta)\vv{z}_k + \beta \vv{x}, {\vv{\lambda}}_k),\label{eq:pp1}\\
\vv{z}_{k+1} & = (1-\beta)\vv{z}_k + \beta \vv{x}_k. \label{eq:pp2}
\end{align}
Then, for $k$ sufficiently large ($k \ge \bar k)$ we have that
$L(\vv{z}_{k+1},\vv{\lambda}_k) - q(\vv{\lambda}_k) \le L(\vv{z}_{k+1},{\vv{\lambda}}_k) - f^* \le 2 \epsilon$
provided  $\beta$ is sufficiently small, \emph{i.e.,} ${0 < } \beta \le (1-\gamma) \min \{ \epsilon/(\bar \mu_L \bar x_D^2),1\}$ where $\bar x_D := 2 \max_{\vv{x} \in D} \| \vv{x} \|_2$, $\bar \mu_L = \mu_f + \bar{\lambda} \vv{1}^T \vv{\mu}_{\vv{g}}$. 
}}  
\end{lemma}
\begin{IEEEproof}
Observe that since
$| L(\vv{z},{\vv{\lambda}}_{k+1})-L(\vv{z},{\vv{\lambda}}_k) | 
 = | ( {\vv{\lambda}}_{k+1}- {\vv{\lambda}}_k)^T\vv{g}(\vv{z}) | 
 \le \|  {\vv{\lambda}}_{k+1}- {\vv{\lambda}}_k\|_2 \| \vv{g}(\vv{z}) \|_2 
 \le \| {\vv{\lambda}}_{k+1}- {\vv{\lambda}}_k\|_2  m{{\bar g}} 
\le  \gamma_1 \gamma \beta \epsilon$
and $L(\cdot,{\vv{\lambda}}_k)$ has uniformly bounded curvature by Theorem \ref{th:descent} we have that for $k$ sufficiently large ($k \ge \bar k$) then $L(\vv{z}_{k+1}, {\vv{\lambda}}_k) -q({\vv{\lambda}}_k) \le 2 \epsilon $ where $q(\vv{\lambda}) := \min_{\vv{z} \in C} L(\vv{z},\vv{\lambda})$. Further, since $q({\vv{\lambda}}) \le q(\vv{\lambda}^*) \le f^*$ for all ${\vv{\lambda}} \succeq \vv{0}$ it follows that $L(\vv{z}_{k+1},{\vv{\lambda}}_k) -f^* \le 2 \epsilon$ for $k \ge \bar k$.
\end{IEEEproof}

\begin{lemma}[Lagrangian of Averages]\label{th:subgradient} 
{\color{black}{Consider optimisation problem $P$ and update
$\vv{\lambda}_{k+1} = [\vv{\lambda}_k + \alpha \vv{g}(\vv{z}_{k+1})]^{[0 , \bar \lambda]}$ where $\alpha > 0$ and $\{\vv{z}_k\}$ is a sequence of points from $C$ such that $L(\vv{z}_{k+1}, \vv{\lambda}_k) - q(\vv{\lambda}_k) \le 2 \epsilon$ for all $k=1,2,\dots$. Let $\vi{\lambda}{j}_1 \in [0, \bar \lambda]$ where $\bar \lambda \ge \vi{\lambda^*}{j}$, $j=1,\dots,m$. Then,
\begin{align}
| L(\vv{z}^\diamond_{k},\vv{\lambda}^\diamond_k) -  f^* |
 \le 2 \epsilon + \frac{\alpha}{2}m {\bar g}^2 + \frac{m \bar \lambda^2}{\alpha k }  \label{eq:lagavg}
\end{align}
where $\vv{z}^\diamond_k := \frac{1}{k} \sum_{i=1}^k \vv{z}_{i+1}$, $\vv{\lambda}^\diamond_k := \frac{1}{k} \sum_{i=1}^k \vv{\lambda}_i$ and ${{\bar g}}:= \max_{\vv{z} \in C} \|\vv{g}(\vv{z})\|_\infty$.}}
\end{lemma}
\begin{IEEEproof}
See Appendix.
\end{IEEEproof}

{\color{black}{
Note that by selecting $\alpha$ sufficiently small in Lemma \ref{th:subgradient} we can obtain a sequence $\{ \vv{\lambda}_k \}$ that changes sufficiently slowly so to satisfy the conditions of Lemma \ref{th:lup}. Further, by Lemma \ref{th:lup} we can construct a sequence of primal variables that satisfy the conditions of Lemma \ref{th:subgradient} for $k \ge \bar k$ and it then follows that (\ref{eq:lagavg}) is satisfied. 
 

Lemma \ref{th:subgradient} requires that $ \vi{\lambda^*}{j} \le \bar \lambda$ for all $j=1,\dots,m$, so it naturally arises the question as to when $\vi{\lambda^*}{j}$ (and so $\bar \lambda$) is bounded. This is clarified in the next lemma, which corresponds to Lemma 1 in \cite{nedicprimal}.
}}

\begin{lemma}[Bounded Multipliers]\label{th:setq} Let $Q_{\delta} := \{\vv{\lambda} \succeq \vv{0}: q(\vv{\lambda}) \ge q(\vv{\lambda}^*) - \delta \}$ with $\delta \ge 0$ and let the Slater condition hold, \emph{i.e.,} there exists a vector $\bar{\vv{z}} \in C$ such that $\vv{g}(\bar{\vv{z}}) \prec \vv{0}$. Then, for every $\vv{\lambda} \in Q_\delta $ we have that 
\begin{align}
\| \vv{\lambda} \|_2 \le \frac{1}{\upsilon} (f(\bar{\vv{z}}) - q(\vv{\lambda}^*) + \delta) \label{eq:setqbound}
\end{align}where $\upsilon := \min_{j \in \{1,\dots,m\}} - \vi{g}{j}(\bar{\vv{z}})$. \end{lemma}
\begin{IEEEproof}
First of all recall that since the Slater condition holds we have strong duality, \emph{i.e.,} $q(\vv{\lambda}^*) = f^*$, and $f^*$ is finite by Proposition 2.1.1. in \cite{convexanalysis}. Now observe that  when 
$\vv{\lambda} \in Q_\delta$ then $q(\vv{\lambda}^*) -\delta \le q(\vv{\lambda}) = \min_{\vv{z} \in C} L(\vv{z}, \vv{\lambda}) \le f(\bar{\vv{z}}) + \vv{\lambda}^T \vv{g}(\bar{\vv{z}})$, and rearranging terms we obtain $-\vv{\lambda}^T \vv{g}(\bar{\vv{z}}) = - \sum_{j=1}^m \vi{\lambda}{j} \vi{g}{j}(\bar{\vv{z}})\le f(\bar{\vv{z}}) - q(\vv{\lambda}^*) + \delta$. Next,  since $\vv{\lambda}\succeq \vv{0}$ and $-\vi{g}{j}(\bar{\vv{z}}) > 0$ for all  $j=1,\dots,m$, let $\upsilon:=\min_{j \in \{1,\dots,m\}} - \vi{g}{j}(\bar{\vv{z}}) $ and see that $\upsilon  \sum_{j=1}^m \vi{\lambda}{j} \le f(\bar{\vv{z}}) - q(\vv{\lambda}^*) + \delta$. Finally, dividing by $\upsilon$ and using the fact that $\| \vv{\lambda} \|_2 \le  \sum_{j=1}^m \vi{\lambda}{j} $ the stated result follows.
\end{IEEEproof}

From Lemma \ref{th:setq} we   have that it is sufficient for $C_0$ to have non-empty relative interior in order for $Q_\delta$ to be a bounded subset in $\mathbb R^m$, and since by definition $\vv{\lambda}^* \in Q_\delta$ then $\vv{\lambda}^*$ is bounded. The bound obtained in Lemma \ref{th:setq} depends on $q(\vv{\lambda}^*)=f^*$, which is usually not known. Nevertheless, we can obtain a looser bound if we use the fact that $-q(\vv{\lambda}^*) \le - q(\vv{\lambda})$ for all $\vv{\lambda} \succeq \vv{0}$. That is, for every $\vv{\lambda} \in Q_\delta$ we have that
\begin{align*}
\| \vv{\lambda} \|_2 \le \frac{1}{\upsilon} (f(\bar{\vv{z}}) - q(\vv{\lambda}_0) + \delta),
\end{align*}where $\vv{\lambda}_0$ is an arbitrary vector in $\mathbb R^m_+$.

{\color{black}{
That is, when the Slater condition is satisfied the upper and lower bounds in (\ref{eq:lagavg}) are finite and can be made arbitrarily small as $k\to \infty$ by selecting the step size $\alpha$ sufficiently small. Convergence of the average of the Lagrangians does not, of course, guarantee that $f(\vv{z}^\diamond_k) \to f^*$ unless we also have complementary slackness, \emph{i.e.,} $ (\vv{\lambda}^\diamond_k)^T \vv{g}(\vv{z}^\diamond_k) \to 0$. Next we present the following lemma, which is a generalisation of Lemma 3 in \cite{nedicprimal}.
}}

\begin{lemma}[Complementary Slackness and Feasibility]\label{th:slackness}
Let the Slater condition hold and suppose $\{ \vv{z}_k \}$ is a sequence of points in $C$ and $\{\tilde{\vv{\lambda}}_k\}$ a sequence of points in $\mathbb{R}^m$ such that (i) $L(\vv{z}_{k+1}, \tilde{\vv{\lambda}}_k ) -q(\tilde{\vv{\lambda}}_k) \le 2\epsilon$ for all $k$ and (ii) $| \vi{\lambda}{j}_k - \vi{\tilde{\lambda}}{j}_k | \le \alpha \sigma_0$, $j=1,\dots,m$ where $\vv{\lambda}_{k+1} = [\vv{\lambda}_k + \alpha \vv{g}(\vv{z}_{k+1})]^+$, $\epsilon \ge 0$, $\alpha > 0$, $\sigma_0 \ge 0$. Suppose also that $\vi{\lambda}{j}_1  \in [0,\bar \lambda]$ with 
\begin{align*}
\textstyle \bar \lambda \ge \frac{3}{\upsilon}(f(\bar{\vv{z}}) - q(\vv{\lambda}^*)+ \delta)  + \alpha m {\bar g}
\end{align*}where $\delta:= \alpha (m{\bar g}^2/2 + m^2\sigma_0 {\bar g}) +2\epsilon$, ${\bar g}:= \max_{\vv{z} \in C} \| \vv{g}(\vv{z})\|_\infty$, $\bar{\vv{z}}$ a Slater vector and $\upsilon := \min_{j \in \{1,\dots,m\}} - \vi{g}{j}(\bar{\vv{z}})$. Then, $\vi{\lambda}{j}_k \le \bar \lambda$ for all $k=1,2,\dots$, {\color{black}{
\begin{align}
 - \frac{m \bar \lambda^2}{2\alpha k} - \frac{\alpha}{2} m {\bar g}^2 
\le (\vv{\lambda}^\diamond_{k})^T \vv{g}(\vv{z}^\diamond_{k}) 
\le  \frac{m \bar \lambda^2}{\alpha k} \label{eq:slackness}
\end{align}and
\begin{align}
\vi{g}{j}(\vv{z}^\diamond_{k}) \le \frac{\bar \lambda}{\alpha k} \label{eq:feasibility}
\end{align}where $\vv{z}^\diamond_{k} := \frac{1}{k} \sum_{i=1}^{k} \vv{z}_{i+1}$ and $\vv{\lambda}^\diamond_{k} := \frac{1}{k} \sum_{i=1}^{ k} \vv{\lambda}_{i}$. 
}}
\end{lemma}
\begin{IEEEproof}
See Appendix. 
\end{IEEEproof}

Lemma \ref{th:slackness} is expressed in a general form where $\tilde{\vv{\lambda}}$ may be any suitable approximation to the usual Lagrange multiplier.  Evidently, the lemma also applies in the special case where $\vv{\lambda}_k = \tilde{\vv{\lambda}}_k$ in which case $\sigma_0 = 0$. {\color{black}{Note from the lemma as well that the running average $\vv{z}^\diamond_k$ is asymptotically attracted to the feasible region as $k$ increases, \emph{i.e.}, $\lim_{k \to \infty} \vv{g}(\vv{z}^\diamond_k) \preceq \vv{0}$}}

We are now in a position to present one of our main results:
\begin{theorem}[Constrained Optimisation] \label{th:maintheorem}
Consider constrained convex optimisation $P$ and the associated sequence of non-convex optimisations $\{ \tilde{P}_k \}$:
\begin{align}
\vv{x}_k & \in \arg \min_{\vv{x} \in D} \ L((1-\beta)\vv{z}_k + \beta\vv{x}, \tilde{\vv{\lambda}}_k) \label{eq:x4thupdate}\\
\vv{z}_{k+1} & = (1-\beta) \vv{z}_k + \beta \vv{x}_k \label{eq:z4thupdate}\\
\vv{\lambda}_{k+1} & =  [\vv{\lambda}_k  + \alpha \vv{g}(\vv{z}_{k+1})]^{[0, \bar \lambda]}\label{eq:minupdate}
\end{align}Let the Slater condition hold and suppose that $| \vi{\lambda}{j}_k - \vi{\tilde \lambda}{j}_k | \le \alpha \sigma_0$ for all $j=1,\dots,m$, $k\ge 1$ with $\sigma_0 \ge 0$. Further, suppose parameters $\alpha $ and $\beta$ are selected sufficiently small that 
\begin{align}
& 0 < \alpha \le \gamma_1 \gamma \beta \epsilon/(m^2 ({\bar g}^2+ 2 \sigma_0{\bar g}) ) \label{eq:alphaselect} \\
& 0 <\beta\le(1-\gamma)\min\{ \epsilon / (\mu_L\bar{x}^2_D),1\} \label{eq:betaselect}
\end{align} with $\epsilon>0$, $\gamma \in (0,1)$, $\gamma_1 \in (0, 1/2)$,  $\bar x_D := 2 \max_{\vv{x} \in D} \| \vv{x} \|_2$, $\mu_L = \mu_f + \bar \lambda \vv{1}^T \vv{\mu}_{\vv{g}}$ and $\bar \lambda$ as given in Lemma \ref{th:slackness}. Then, for every $\epsilon > 0$ and for $k$ sufficiently large $(k \ge \bar k)$ the sequence of solutions $\{\vv{z}_k\}$ to sequence of optimisations $\{\tilde{P}_k\}$ satisfies: {\color{black}{
 \begin{align}
& - \frac{2 m \bar \lambda^2}{\alpha k} -  \alpha (m{\bar g}^2/2 + m^2\sigma_0 \bar g)  - 2 \epsilon \notag \\
& \quad   \le f(\vv{z}^\diamond_{k}) -  f^* 
 \le 2 \epsilon + \alpha (m{\bar g}^2 + m^2 \sigma_0 \bar g) + \frac{3 m \bar \lambda^2}{2 \alpha k }  \label{eq:mainbound}
 \end{align}where ${\vv{z}}^\diamond_k := \frac{1}{k} \sum_{i=\bar k}^{\bar k + k} {\vv{z}}_{i+1}$, $\tilde{\vv{\lambda}}^\diamond_k := \frac{1}{k} \sum_{i=\bar k}^{\bar k + k} \tilde{\vv{\lambda}}_i$ and ${\bar g} :=  \max_{\vv{z} \in C} \| \vv{g}(\vv{z})\|_\infty$.
 }}
\end{theorem}

\begin{IEEEproof}
{\color{black}{
First of all observe that since $\vi{\lambda}{j}_{k+1} = [\vi{\lambda}{j}_k + \alpha \vi{g}{j}(\vv{z}_{k+1})]^{[0,\bar \lambda]}$ we have that $| \vi{\lambda}{j}_{k+1} - \vi{\lambda}{j}_k | \le \alpha {\bar g}$ for all $k$. Further, since $| \vi{\lambda}{j}_k - \vi{\tilde \lambda}{j}_k  | \le \alpha \sigma_0$ then $
 | \vi{\tilde\lambda}{j}_{k+1} - \vi{\tilde\lambda}{j}_k  | 
 = | \vi{\tilde\lambda}{j}_{k+1} - \vi{\tilde\lambda}{j}_k  + \vi{\lambda}{j}_{k+1} - \vi{\lambda}{j}_{k+1} + \vi{\lambda}{j}_{k} - \vi{\lambda}{j}_k |  
 \le | \vi{\tilde\lambda}{j}_{k+1} - \vi{\lambda}{j}_{k+1}|  +  |\vi{\lambda}{j}_{k+1}  - \vi{\lambda}{j}_{k}| + |\vi{\lambda}{j}_k  - \vi{\tilde\lambda}{j}_k|  
 \le \alpha (2 \sigma_0 + {\bar g})$. That is, 
 \begin{align}
 \| \tilde{\vv{\lambda}}_{k+1} - \tilde{\vv{\lambda}}_k \|_2 \le \alpha m  (2 \sigma_0 + {\bar g}) \qquad k=1,2,\dots
 \end{align} Next, see that since $L(\cdot,\vv{\lambda}_k)$ has uniform bounded curvature and $| L(\vv{z},\tilde{\vv{\lambda}}_{k+1})-L(\vv{z},\tilde{\vv{\lambda}}_k) |  \le \|\tilde{\vv{\lambda}}_{k+1}-\tilde{\vv{\lambda}}_k\|_2 \| \vv{g}(\vv{z}_{k+1}) \|_2 \le \|\tilde{\vv{\lambda}}_{k+1}-\tilde{\vv{\lambda}}_k\|_2 m {{\bar g}} \le \alpha m^2  (2 \sigma_0 {\bar g}+ {\bar g}^2)  \le  \gamma_1 \gamma \beta \epsilon$, it follows by Lemma \ref{th:lup} that for $k$ sufficiently large ($k \ge \bar k$) then 
 $
 L(\vv{z}_{k+1},\tilde{\vv{\lambda}}_k) - q(\tilde{\vv{\lambda}}_k) \le 2 \epsilon
 $
and therefore by Lemma \ref{th:subgradient} 
\begin{align*}
& - \frac{m \bar \lambda^2}{\alpha k} -  \frac{\alpha}{2} m {\bar g}^2  - 2 \epsilon \notag \\
& \quad \qquad  \le L(\vv{z}^\diamond_{k},\tilde{\vv{\lambda}}^\diamond_k) -  f^* 
 \le 2 \epsilon +  \frac{\alpha}{2}m {\bar g}^2 + \frac{m \bar \lambda^2}{\alpha k }. 
\end{align*}
Next, see that since 
$ | L(\vv{z}^\diamond_{k},{\vv{\lambda}}^\diamond_k) -  L(\vv{z}^\diamond_{k}, \tilde{\vv{\lambda}}^\diamond_k) | =  ({\vv{\lambda}}^\diamond_k - \tilde{\vv{\lambda}}^\diamond_k)^T \vv{g}(\vv{z}^\diamond_k)  \le \| \vv{\lambda}^\diamond_k - \tilde{\vv{\lambda}}^\diamond_k \|_2 \| \vv{g}(\vv{z}^\diamond_k) \|_2 \le \alpha m^2 \sigma_0 \bar g$ we have that 
\begin{align*}
& - \frac{m \bar \lambda^2}{\alpha k} -  \alpha (m{\bar g}^2 /2 + m^2\sigma_0 \bar g)  - 2 \epsilon \notag \\
& \quad \quad  \le L(\vv{z}^\diamond_{k},{\vv{\lambda}}^\diamond_k) -  f^* 
 \le 2 \epsilon + \alpha (m{\bar g}^2/2 + m^2 \sigma_0 \bar g) + \frac{m \bar \lambda^2}{\alpha k } .
\end{align*}Finally, by using the complementary slackness bound of Lemma \ref{th:slackness} the stated result follows. 
}}
\end{IEEEproof}

Theorem \ref{th:maintheorem} says that by selecting step size $\alpha$ and smoothing parameter $\beta$ sufficiently small then the average of the solutions to the sequence of non-convex optimisations $\{ \tilde{P}_k \}$ can be made arbitrarily close to the solution of constrained convex optimisation $P$.


\subsubsection{Alternative Update}

Note that, by replacing use of Theorem \ref{th:descent} by Theorem \ref{th:descent1b} in the proof, we can replace update (\ref{eq:x4thupdate}) by its non-convex Frank-Wolfe alternative,
\begin{align}
\vv{x}_k & \in \arg \min_{\vv{x} \in D} \ \partial_{\vv{z}} L(\vv{z}_k, \vv{\lambda}_k)^T \vv{x} \label{eq:x4thupdatefw} \\
&= \arg \min_{\vv{x} \in D} \ (\partial f(\vv{z}_k) +  \vv{\lambda}_k^T\partial \vv{g}(\vv{z}_k))^T \vv{x}. \notag
\end{align}
That is, we have:
\begin{corollary}[Constrained Optimisation Using Frank-Wolfe Update]\label{th:convconvergencefw}
Consider the setup in Theorem \ref{th:maintheorem} but with update (\ref{eq:x4thupdate})  replaced by (\ref{eq:x4thupdatefw}). Then, there exists a finite $\bar k$ such that the bound given in (\ref{eq:mainbound}) holds. 
\end{corollary}


\subsection{Generalised Update} 

Let $C^\prime\subseteq \conv(D)$ be any subset of the convex hull of action set $D$, including the empty set.   Since $\min_{\vv{x} \in C^\prime \cup D} \ L((1-\beta)\vv{z}_k + \beta\vv{x}, \vv{\lambda}_k) \le \min_{\vv{x} \in D} \ L((1-\beta)\vv{z}_k + \beta\vv{x}, \vv{\lambda}_k)$, we can immediately generalise update (\ref{eq:x4thupdate}) to
\begin{align}
\vv{x}_k & \in \underset{\vv{x} \in C^\prime\cup D}{\arg \min} \ L((1-\beta)\vv{z}_k + \beta\vv{x}, \vv{\lambda}_k) \label{eq:gen1}
\end{align}
and Theorem \ref{th:maintheorem} will continue to apply.   Selecting $C^\prime$ equal to the empty set we recover (\ref{eq:x4thupdate}) as a special case.  Selecting $C^\prime=\conv(D)$ we recover the classical convex dual subgradient update as a special case. Update (\ref{eq:gen1}) therefore naturally generalises both the classical convex dual subgradient update and non-convex update  (\ref{eq:x4thupdate}). Hence, we have the following corollary.

\begin{corollary}[Constrained Optimisation Using Unified Update] \label{th:corogenupdate}
Consider the setup in Theorem \ref{th:maintheorem} but with update (\ref{eq:x4thupdate})  replaced by (\ref{eq:gen1}). Then, there exists a finite $\bar k$ such that the bound given in (\ref{eq:mainbound}) holds. 
\end{corollary}


\section{Using Queues As Approximate Multipliers}
\label{sec:perturbedmultiplier}

In Theorem \ref{th:maintheorem} the only requirement on the sequence of approximate multipliers $\{\tilde{\vv{\lambda}}_k\}$ is that it remains close to the sequence of Lagrange multipliers $\{\vv{\lambda}_k\}$ generated by a dual subgradient update in the sense that $| \vi{\lambda}{j}_k - \vi{\tilde{\lambda}}{j}_k | \le \alpha \sigma_0$ for all $k$.  In this section we consider the special case where sequence $\{\tilde{\vv{\lambda}}_k\}$ additionally satisfies the following,
\begin{align}
\tilde{\vv{\lambda}}_{k+1} = [\tilde{\vv{\lambda}}_k + \tilde{\vv{\delta}}_k]^{[0, \bar \lambda]} \label{eq:perturbedupdateorig}
\end{align}
with $\tilde{\vv{\delta}}_k\in\mathbb{R}^m$ and $\tilde{\vv{\lambda}}_1 = \vv{\lambda}_1$.

We begin by recalling the following lemma, which is a direct result of \cite[Proposition 3.1.2]{meyn2008control}.
\begin{lemma} \label{th:sequences}
Consider sequences $\{ \lambda_k\}$, $\{ \tilde \lambda_k\}$ in $\mathbb R$ given by updates $\lambda_{k+1} = [\lambda_k +  \delta_k]^{[0, \bar \lambda]}$, $ \tilde \lambda_{k+1} = [ \tilde \lambda_k +  \tilde \delta_k]^{[0, \bar \lambda]}$ where $\delta, \tilde \delta \in \mathbb R$. Suppose $\lambda_1  = \tilde \lambda_1$ and $| \sum_{i=1}^k \delta_i - \tilde \delta_i | \le \epsilon$ for all $k$. Then, 
\begin{align*}
| \lambda_{k}  - \tilde{\lambda}_{k} | \le  2 \epsilon \qquad k=1,2,\dots
\end{align*}
\end{lemma}
\begin{IEEEproof}
See Appendix.
\end{IEEEproof}
Applying Lemma \ref{th:sequences} to our present context it follows that $| \vi{\lambda}{j}_k - \vi{\tilde{\lambda}}{j}_k | \le \alpha \sigma_0$ for all $k$ (and so Theorem \ref{th:maintheorem} holds) for every sequence $\{ \tilde{\vv{\delta}}_k\}$ such that $| \sum_{i=1}^k \alpha \vi{g}{j}(\vv{z}_{i}) - \vi{\tilde{\delta}}{j}_i | \le \alpha \sigma_0$ for all $k$.

Of particular interest is the special case of optimisation $P$ where the constraints are linear.   {That is, ${\vi{g}{j}}(\vv{z})=\vvi{a}{j} \vv{z}-\vi{b}{j}$ where $(\vvi{a}{j})^T \in\mathbb{R}^{n}$ and $\vi{b}{j}\in\mathbb{R}$, $j=1,\dots, m$. Gathering vectors $\vvi{a}{j}$ together as the rows of matrix $\vv{A}\in\mathbb{R}^{m\times n}$} and collecting additive terms $\vi{b}{j}$ into vector $\vv{b}\in\mathbb{R}^m$, the linear constraints can then be written as $\vv{A}\vv{z}\preceq \vv{b}$. Therefore, the dual subgradient Lagrange multiplier update in the sequence of optimisations $\{\tilde P _k\}$ is given by
\begin{align}
{\vv{\lambda}}_{k+1}  & = [ {\vv{\lambda}}_k + \alpha (\vv{A}{\vv{z}}_{k+1} - \vv{b})]^{[0 , \bar \lambda]} \label{eq:lzupdate}
\end{align}
with $\vv{z}_{k+1} = (1-\beta)\vv{z}_k + \beta \vv{x}_k$, $\vv{x}_k \in D$. 
Now suppose that in (\ref{eq:perturbedupdateorig}) we select $\tilde{\vv{\delta}}_k=\alpha (\vv{A}\vv{x}_{k} - \vv{b}_k)$ where $\{ \vv{b}_k \}$ is a sequence of points in $\mathbb R^m$.  Then,
\begin{align}
\tilde{\vv{\lambda}}_{k+1}  & =  [\tilde{\vv{\lambda}}_k + \alpha (\vv{A}\vv{x}_{k} - \vv{b}_k)]^{[0 , \bar \lambda]}\label{eq:lxupdate}
\end{align}
with $\tilde{\vv{\lambda}}_1 = \vv{\lambda}_1$. 

Observe that in (\ref{eq:lxupdate}) we have replaced the continuous-valued quantity $\vv{z}_{k}$ with the discrete-valued quantity $\vv{x}_{k}$.  We have also replaced the constant $\vv{b}$ with the time-varying quantity $\vv{b}_k$.  Further, letting
$\vv{Q} : = \tilde{\vv{\lambda}} / \alpha$ then (\ref{eq:lxupdate}) can be rewritten equivalently as
\begin{align}
{\vv{Q}}_{k+1}  & = [ {\vv{Q}}_k + \vv{A}{\vv{x}}_{k} - \vv{b}_k]^{[0, \bar \lambda / \alpha]} \label{eq:queueupper}
\end{align}
which is a discrete queue length update with increment $\vv{A}{\vv{x}}_{k} - \vv{b}_k$.   The approximate multipliers $\tilde{\vv{\lambda}}$ are therefore scaled discrete queue occupancies.

Using Lemma \ref{th:sequences}  it follows immediately that Theorem \ref{th:maintheorem} holds provided
\begin{align}
\textstyle |\sum_{i=1}^k \vvi{a}{j}(\vv{z}_i-\vv{x}_i) +  (\vi{b}{j}_i-\vi{b}{j}) |  \le \alpha \sigma_0\label{eq:cond}
\end{align}
Since update $\vv{z}_{k+1} = (1-\beta)\vv{z}_k + \beta \vv{x}_k$ yields a running average of $\{\vv{x}_k\}$ we might expect that sequences $\{\vv{z}_k\}$ and $\{\vv{x}_k\}$ are always close and so uniform boundedness of $| \sum_{i=1}^k (\vi{b}{j}_i-\vi{b}{j})|$ is sufficient to ensure that (\ref{eq:cond}) is satisfied.   This is indeed the case, as established by the following theorem.



\begin{theorem}[Queues as approximate multipliers] \label{th:auxiliarymultiplier}
Consider updates (\ref{eq:lzupdate}) and (\ref{eq:lxupdate}) where $\{ \vv{x}_k \}$ is  an arbitrary sequence of points in $D$, $\vv{z}_{k+1}  = (1-\beta) \vv{z}_k + \beta \vv{x}_k$, $\beta \in (0,1)$, $\vv{z}_1 \in C:=\conv{(D)}$. Further, suppose that $\{\vv{b}_k\}$ is a sequence of points in $\mathbb R^m$ such that $| \sum_{i=1}^k (\vi{b}{j}_i - \vi{b}{j}) | \le \sigma_2$  for all $j=1,\dots,m$, $k=1,2,\dots$. Then,
\begin{align*}
\| \tilde{\vv{\lambda}}_{k} - \vv{{\lambda}}_{k}\|_2 \le 2m \alpha (\sigma_1 / \beta + \sigma_2) ,\quad k=1,2,\dots
\end{align*}
where ${\sigma_1} :=  2 \max_{\vv{z} \in C} \| \vv{A} \vv{z} \|_\infty$. 
\end{theorem}
\begin{IEEEproof}
See Appendix.
\end{IEEEproof}
Observe that the difference between $\vv{\lambda}_k$ and $\tilde{\vv{\lambda}}_k$ can be made arbitrarily small by selecting $\alpha$ small enough.  The requirement that $| \sum_{i=1}^k (\vi{b}{j}_i - \vi{b}{j}) | \le \sigma_2$ is satisfied when sequence $\{ \vi{b}{j}_k \}$ converges sufficiently fast to $\vi{b}{j}$ (dividing both sides by $k$, the requirement is that $|\frac{1}{k} \sum_{i=1}^k \vi{b}{j}_i - \vi{b}{j} | \le \sigma_2 / k$).   

\begin{figure}
\centering
\includegraphics[width=\columnwidth]{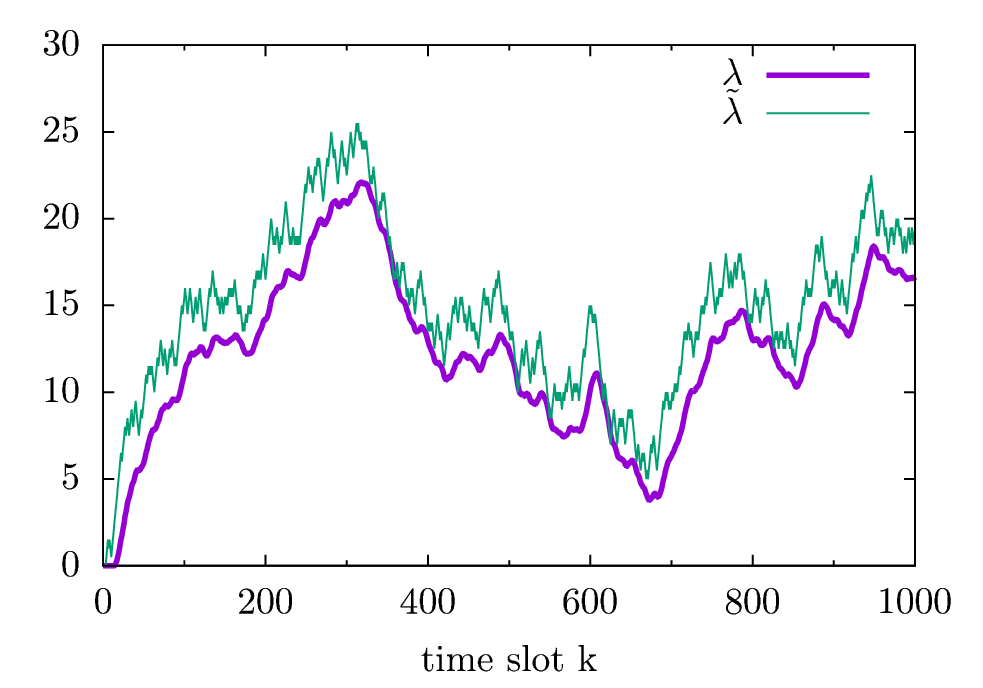}
\caption{Example realisations of $\tilde \lambda_k$ (thin line) and ${\lambda}_k$ (thicker line) given by updates (\ref{eq:lzupdate}) and (\ref{eq:lxupdate}).} \label{fig:q}
\end{figure}

In the special case when $\vi{b}{j}_k =\vi{b}{j}$ then Theorem \ref{th:auxiliarymultiplier} is trivially satisfied.  This is illustrated in Figure \ref{fig:q}, which plots $\vv{\lambda}_k$ and $\tilde{\vv{\lambda}}_k$ for a simple example where $A = 1$, $b_k = b=0.5$, $\alpha=1$, $\beta=0.1$ and sequence $\{x_k\}$ takes values independently and uniformly at random from set $\{0,1\}$.  It can be seen that the distance between $\vv{\lambda}_k$ and $\tilde{\vv{\lambda}}_k$ remains uniformly bounded over time. 

In summary, we have arrived at the following corollary to Theorem \ref{th:maintheorem}.

\begin{corollary}[Constrained Optimisation Using Approximate Multipliers] \label{th:lincoro}
Consider the setup of Theorem \ref{th:maintheorem}, suppose the constraints are linear $\vv{A}\vv{z} - \vv{b}\preceq \vv{0}$ and $\tilde{\vv{\lambda}}_{k+1} = [\tilde{\vv{\lambda}}_k + \alpha(\vv{A} \vv{x}_k - \vv{b}_k)]^{[0, \bar \lambda]}$, $\vv{b}_k\in\mathbb{R}^m$. Suppose $| \frac{1}{k} \sum_{i=1}^k \vi{b}{j}_i - \vi{b}{j}| \le \sigma_2 / k$ for all $j$ and $k$. Then, the bound (\ref{eq:mainbound}) in Theorem \ref{th:maintheorem} holds with $\sigma_0 = 2(\sigma_1/\beta + \sigma_2)$  where ${\sigma_1} :=  2 \max_{\vv{z} \in C} \| \vv{A} \vv{z} \|_\infty$. 
\end{corollary}

\subsection{Weaker Condition for Loose Constraints}
{
Suppose constraint $j$ is loose at an optimum, \emph{i.e.}, $\vi{g}{j}(\vv{z}^*)<0$ for $\vv{z}^*\in C^*$.   Then by complementary slackness the associated Lagrange multiplier must be zero, \emph{i.e.}, $(\vi{\lambda}{j})^*=0$, and we can select $\vi{\lambda}{j}_k = (\vi{\lambda}{j})^* = 0$ for all $k$. Since $\vi{\tilde\lambda}{j}_k$ is non-negative, to apply Theorem \ref{th:maintheorem} it is enough that $\vi{\tilde\lambda}{j}_k \le \alpha\sigma_0$ for $k=1,2,\dots$.  Assuming, for simplicity, that $\vi{\tilde \lambda}{j}_1=0$, from the proof of Lemma \ref{th:sequences} we have $\vi{\tilde\lambda}{j}_k = [ \max_{1\le l \le k-1}\sum_{i=l}^{k-1} \alpha(\vvi{a}{j}\vv{x}_i -  \vi{b}{j}_i )]^+$
and so a sufficient condition for $\vi{\tilde\lambda}{j}_k \le \alpha\sigma_0$ is that 
$\max_{1\le l \le {k-1}}\sum_{i=l}^{k-1} (\vvi{a}{j}\vv{x}_i -\vi{b}{j}) -  (\vi{b}{j}_i -\vi{b}{j})\le \sigma_0$ for all $k$.  
The advantage of this condition is that $-\sum_{i=l}^{k-1} (\vi{b}{j}_i -\vi{b}{j})$ is now not bounded below and so a wider class of sequences $\{\vi{b}{j}_i\}$ is potentially admissible.  The disadvantage is that to exploit this we need to know in advance that constraint $j$ is loose at the optimum.}

\subsection{Queue Stability}
{{Recall that by Lemma \ref{th:slackness} sequence $\{\vv{\lambda}_k\}$ in Theorem \ref{th:maintheorem} (and respective corollaries of the theorem)}} is bounded for all $k \ge \bar k$. Therefore, since $\| \tilde{\vv{\lambda}}_{k} - \vv{{\lambda}}_{k}\|_2$ is uniformly bounded it follows that $\{\tilde{\vv{\lambda}}_k\}$ is also bounded and therefore the associated discrete queue is stable (although the occupancy $\vv{Q}$ of the discrete queue scales with $1/\alpha$ since $\vv{Q} = \tilde{\vv{\lambda}} / \alpha$).   Note that we have arrived to this queue stability result purely from a convex optimisation analysis and without using any Foster-Lyapunov argument.


\subsection{Optimal Actions Depend Only on Queue Occupancy}
In network resource allocation problems where the linear constraints can be identified with link queues we can use the scaled queue occupancies directly in the optimisation. That is, 
\begin{align}
\vv{x}_k 
& \in \arg \min_{\vv{x} \in D} L((1-\beta)\vv{z}_k + \beta \vv{x}, \alpha {\vv{Q}}_k) \label{eq:orig} \\
& =  \arg \min_{\vv{x} \in D} f((1-\beta)\vv{z}_k + \beta \vv{x})  + \alpha \beta {\vv{Q}}_k^T \vv{A}\vv{x} \label{eq:simple}
\end{align}where update (\ref{eq:simple}) is obtained from (\ref{eq:orig}) by retaining only the parts of $L((1-\beta)\vv{z}_k + \beta\vv{x}, \tilde{\vv{\lambda}}_k)$ which depend on $\vv{x}$ \emph{i.e.,} dropping constant terms which do not change the solution to the optimisation. We could also consider Corollary \ref{th:convconvergencefw} and so have a Frank-Wolfe like update: 
\begin{align}
\vv{x}_k 
& \in \arg \min_{\vv{x} \in D} \partial_{\vv{z}} L(\vv{z}_k, \alpha {\vv{Q}}_k)^T \vv{x} \label{eq:origfw} \\
& =  \arg \min_{\vv{x} \in D} \partial f(\vv{z}_k)^T \vv{x}  + \alpha {\vv{Q}}_k^T \vv{A}\vv{x} \label{eq:simplefw}
\end{align}
Importantly, note that neither (\ref{eq:simple}) nor (\ref{eq:simplefw}) involve $\vv{b}$ or $\vv{b}_k$. Therefore, we can generate a sequence of discrete actions by simply looking at the queue occupancies at each time slot.

\section{Stochastic Optimisation}

The analysis in the preceeding Section \ref{sec:convc} is for deterministic optimisation problems.   However, it can be readily extended to a class of stochastic optimisations.   

%
%
%
%

\subsection{Stochastic Approximate Multipliers}

\subsubsection{Linear Constraints}
Of particular interest, in view of the equivalence which has been established between approximate multipliers and queues, is accommodating stochastic queue arrivals.  

Let $\{\vv{B}_k\}$ be a stochastic process with realisations of $\vv{B}_k$ taking values in $\mathbb{R}^m$ and with mean $\vv{b}$. Let $p_K:=\Prob(\max_{ k\in\{1,2,\dots,K\}}\|\sum_{i=1}^k (\vv{B}_i -\vv{b})\|_{\infty} \le {\sigma_2})$. Let  $\{\vv{b}_i\}_{i=1}^K$ denote a  realisation of length $K$. Fraction $p_K$ of these realisations satisfy $\|\sum_{i=1}^k (\vv{b}_i-\vv{b})\|_{\infty} \le {\sigma_2}$ for all $k=1,2,\dots,K$.  When this fraction is asymptotically lower bounded $\liminf_{K\rightarrow\infty}p_K\ge p$, then fraction $p$ of realisations satisfy the conditions of Theorem \ref{th:auxiliarymultiplier}.  We therefore have the following corollary (which is a stochastic version of Corollary \ref{th:lincoro}) to Theorem \ref{th:maintheorem}.

\begin{corollary}
Consider the setup of Theorem \ref{th:maintheorem}, suppose the constraints are linear $\vv{A}\vv{z} - \vv{b}\preceq \vv{0}$ and $\tilde{\vv{\lambda}}_{k+1} = [\tilde{\vv{\lambda}}_k + \alpha(\vv{A} \vv{x}_k - \vv{b}_k)]^{[0, \bar \lambda]}$, $\vv{b}_k\in\mathbb{R}^m$, $\tilde{\vv{\lambda}}_1 = \vv{\lambda}_1$.  Let sequence $\{\vv{b}_k\}$ be a realisation of a stochastic process $\{\vv{B}_k\}$ and $p_K:=\Prob(\max_{ k\in\{1,2,\dots,K\}}\|\sum_{i=1}^k(\vv{B}_i-\vv{b})\|_{\infty} \le {\sigma_2})$.  Suppose that this probability is asymptotically lower bounded $\liminf_{K\rightarrow\infty}p_K\ge p$.   Then, with probability at least $p$ the bound (\ref{eq:mainbound}) in Theorem \ref{th:maintheorem} holds with $\sigma_0 = 2(\sigma_1/\beta + \sigma_2)$, ${\sigma_1} :=  2 \max_{\vv{z} \in C} \| \vv{A} \vv{z} \|_\infty$.
\end{corollary}

Note that there is no requirement for stochastic process $\{\vv{B}_k\}$ to be i.i.d. or for any of its properties, other than that feasible set $\vv{A}\vv{z}\preceq\vv{b}$ has non-empty relative interior, to be known in advance in order to construct solution sequence $\{\tilde{P}_k\}$. 

\subsubsection{Nonlinear Constraints}
We can further generalise the latter corollary to consider non-linear stochastic constraints: 

\begin{corollary} \label{th:seq}
Consider the setup in Theorem \ref{th:maintheorem} with update $\tilde{\vv{\lambda}}_{k+1} = [ \tilde{\vv{\lambda}}_k + \alpha(\vv{g}(\vv{z}_{k+1}) - \vv{b}_k)]^{[0, \bar \lambda]}$, $\vv{b}_k \in \mathbb R^m$, $\tilde{\vv{\lambda}}_1 = \vv{\lambda}_1$. Let sequence $\{\vv{b}_k\}$ is a realisation of a stochastic process $\{\vv{B}_k\}$ with $\vv{0}$ mean  and $p_K:=\Prob(\max_{ k\in\{1,2,\dots,K\}}\|\sum_{i=1}^k(\vv{B}_i-\vv{b})\|_{\infty} \le {\sigma_2})$.  Suppose that this probability is asymptotically lower bounded $\liminf_{K\rightarrow\infty}p_K\ge p$.   Then, with probability at least $p$ the bound (\ref{eq:mainbound}) in Theorem \ref{th:maintheorem} holds with $\sigma_0 = 2\sigma_2$.
\end{corollary}




\subsection{Stochastic Actions}

Suppose that when at time $k$ we select action $\vv{x}_k\in D$, the action actually applied is a realisation of random variable $\vv{Y}_k$ that also takes values in $D$; this is for simplicity, the extension to random action sets different from $D$ is straightforward.  For example, we may select $x_k=1$ (which might correspond to transmitting a packet) but with some probability actually apply $y_k=0$ (which might correspond to a transmission failure/packet loss).  Let $p_{\vv{x}\vv{y}}:=\Prob(\vv{Y}_k=\vv{y} | \vv{x}_k=\vv{x})$, $\vv{x}$, $\vv{y}\in D$ and we assume that this probability distribution is time-invariant \emph{i.e.,} does not depend on $k$; again, this can be relaxed in the obvious manner.   

Namely, assume that the probabilities $p_{\vv{x}\vv{y}}$, $\vv{x}$, $\vv{y}\in D$ are known.   Then $\bar{\vv{y}}(\vv{x}) := \mathbb{E}[\vv{Y}_k|\vv{x}_k=\vv{x}] = \sum_{\vv{y}\in D} \vv{y} p_{\vv{x}\vv{y}}$ can be calculated.   The above analysis now carries over unchanged provided we modify the non-convex optimisation from $\min_{\vv{x} \in D} L((1-\beta)\vv{z}_k + \beta \vv{x},{\vv{\lambda}}_k)$ to $ \min_{\vv{x} \in D} L((1-\beta)\vv{z}_k + \beta \bar{\vv{y}}(\vv{x}),{\vv{\lambda}}_k)$ and everywhere replace $\vv{x}_k$ by $\bar{\vv{y}}(\vv{x}_k)$.  That is, we simply change variables to $\bar{\vv{y}}$.  Note that this relies upon the mapping from $\vv{x}$ to $\bar{\vv{y}}$ being known.  If this is not the case, then we are entering the realm of stochastic decision problems and we leave this to future work.

\section{Max-Weight Revisited}

\subsection{Discussion}
Recall the formulation of a queueing network in Section \ref{sec:introduction}, where matrix $\vv{A}$ defines the queue interconnection, with $i$'th row having a $-1$ at the $i$'th entry, $1$ at entries corresponding to queues from which packets are sent to queue $i$, and $0$ entries elsewhere.  Hence, the queue occupancy evolves as $\vv{Q}_{k+1}=[\vv{Q}_k+\vv{A}\vv{x}_k+\vv{b}_k]^{[0,\bar \lambda / \alpha]}$.  As shown in Section \ref{sec:perturbedmultiplier}  updates $\vv{x}_k \in \arg\min_{\vv{x}\in D} \partial f(\vv{z}_k)^T \vv{x} + \alpha\vv{Q}^T_k\vv{A}\vv{x}$, $\vv{z}_{k+1}=(1-\beta)\vv{z}_k+\beta\vv{x}_k$ leads to $\vv{z}_k$ converging to a ball around the solution to the following convex optimisation,
\begin{align*}
\underset{\vv{z}\in C}{\text{minimise}} \qquad &  f(\vv{z})\\
\underset{}{\text{subject to}} \qquad & \vv{A}\vv{z}+\vv{b}\preceq 0
\end{align*}
where $C=\conv(D)$, $\{\vv{b}_k\}$ is any sequence such that $\lim_{k\rightarrow\infty}\frac{1}{k}\sum_{i=1}^k\vv{b}_i = \vv{b}$ and $|(\frac{1}{k}\sum_{i=1}^k\vi{b}{j}_i)-\vi{b}{j}|\le {\sigma_2}/k$, $j=1,\dots,m$ for some finite $\sigma_2>0$.  

Observe that this update is identical to the greedy primal-dual max-weight schedule once we identify utility function $U(\cdot)$ with $-f(\cdot)$.   However, we have arrived at this from a purely convex optimisation viewpoint and by elementary arguments, without recourse to more sophisticated Lyapunov drift, stochastic queueing theory \emph{etc}. Further, our analysis immediately generalises the max-weight analysis to allow arbitrary linear constraints rather than just the specific constraints associated with a queueing network, and beyond this to convex nonlinear constraints with bounded curvature.  

In our analysis, the key role played by bounded curvature in non-convex descent is brought to the fore.  This property is of course present in existing max-weight results, in the form of a requirement for continuous differentiability of the utility function, but insight into the fundamental nature of this requirement had been lacking.  One immediate benefit is the resulting observation that any non-convex update with suitable descent properties can be used, and strong connections are established with the wealth of convex descent methods.  For example, by Theorem \ref{th:maintheorem} we can replace update $\vv{x}_k \in \arg\min_{\vv{x}\in D}  (\partial f(\vv{z}_k) + \vv{A}\vv{Q}_k)^T\vv{x}$ (which is now seen to be a variant of the classical Frank-Wolfe update) with the  direct Lagrangian update $\vv{x}_k \in \arg\min_{\vv{x}\in D} f(\vv{z}_k + \beta(\vv{x}-\vv{z}_k)) +  \beta \vv{Q}^T_k\vv{A}\vv{x}$ to obtain a new class of non-convex algorithms.

\section{Numerical Examples}

\subsection {Example: Convergence and Bounds in Theorem \ref{th:maintheorem}} \label{ex:bt3ex}

Consider the convex optimisation problem $\min_{\vv{z}\in C}\sum_{i=1}^n \exp(\vv{V} \vv{z})$ s.t. $\vv{b} \preceq \vv{z}$ where $\vv{V}:= \mathrm{diag}(1,\dots,n)$, $C:= \conv{(D)}$, $D:= \{\vv{x} \in \mathbb R^n : \vi{x}{j} \in \{0,s\}, j=1,\dots,n \}$, $s > 0$ and $\vv{b} = ({s}/{\sum_{i=1}^n 2i}) [1,\dots,n]^T$.  Observe that  the Slater condition holds.  Consider the following sequence of non-convex optimisations for $k=1,2,\dots$,
\begin{align}
\vv{x}_k & \in \arg \min_{\vv{x} \in D} \sum_{i=1}^n \exp(\vv{V}((1-\beta)\vv{z}_k + \beta\vv{x})) \notag \\
& \qquad \qquad \qquad + \tilde{\vv{\lambda}}^T_k (\vv{b} - ((1-\beta)\vv{z}_k + \beta\vv{x}))  \label{eq:upex1}\\
\vv{z}_{k+1} & = (1-\beta)\vv{z}_k + \beta \vv{x}_k \label{eq:upex2} \\
\vv{\lambda}_{k+1} & = [ \vv{\lambda}_k + \alpha (\vv{b} - \vv{z}_{k+1}) ]^{[0, \bar \lambda]} \label{eq:upex3}
\end{align}
with $\vv{z}_1 = s\vv{1}$, $\vi{\lambda}{j}_1, \vi{\tilde{\lambda}}{j}_1 = 0$, $j=1,\dots,m$ and parameters $\alpha$ and $\beta$ are selected as indicated in (\ref{eq:alphaselect}) and (\ref{eq:betaselect}), with parameters $n=3$, $s = 1/\sqrt{n \bar{\mu}_L}$, $\bar \mu_L = 0.6$, $\gamma = 0.5$, $\bar \lambda = 0.7$, ${\bar g} = 0.6211$ and  $\bar x_D = s \sqrt{n}$. 

\begin{figure}
\centering
\includegraphics[width=\columnwidth]{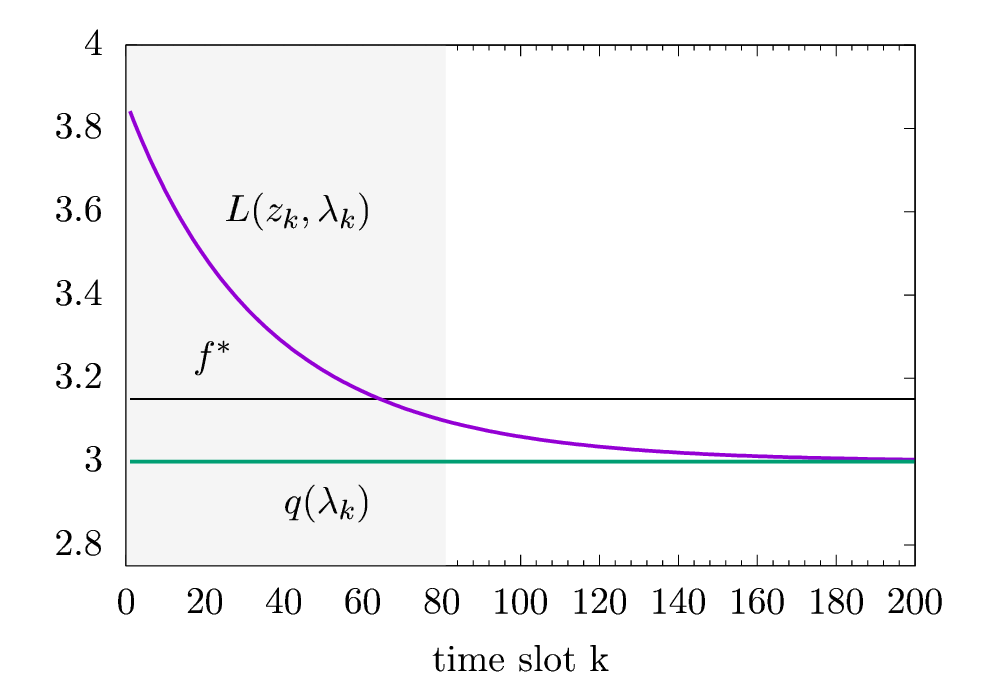}
\caption{Illustrating the convergence of $L(\vv{z}_{k+1} , \vv{\lambda}_k)$ to a ball around $q(\vv{\lambda}_k)$ for $\epsilon = 0.05$ and $\sigma_0 = 0$.  Shaded area ($k < 81$) indicates that $L(\vv{z}_{k+1}, \vv{\lambda}_{k}) - q(\vv{\lambda}_k) > 2 \epsilon$.}
\label{fig:2econv}
\end{figure}

\subsubsection{Convergence into $2\epsilon$-ball in finite time}
To begin with, suppose ${\tilde{\lambda}} = {\lambda}$.  For $\epsilon = 0.05$ (so $\alpha = 7.29 \cdot 10^{-5}$) Figure \ref{fig:2econv} plots the convergence of $L(\vv{z}_{k+1} , \vv{\lambda}_k)$ into an $2\epsilon$-ball around $f^*$.   It can be seen that this convergence occurs within finite time, $\bar k = 81$ and that $L(\vv{z}_{k+1} , \vv{\lambda}_k)$ then stays within this ball at times $k\ge  \bar k$.

\subsubsection{Upper and lower bounds from Theorem \ref{th:maintheorem}}
\begin{figure}
\centering
\includegraphics[width=\columnwidth]{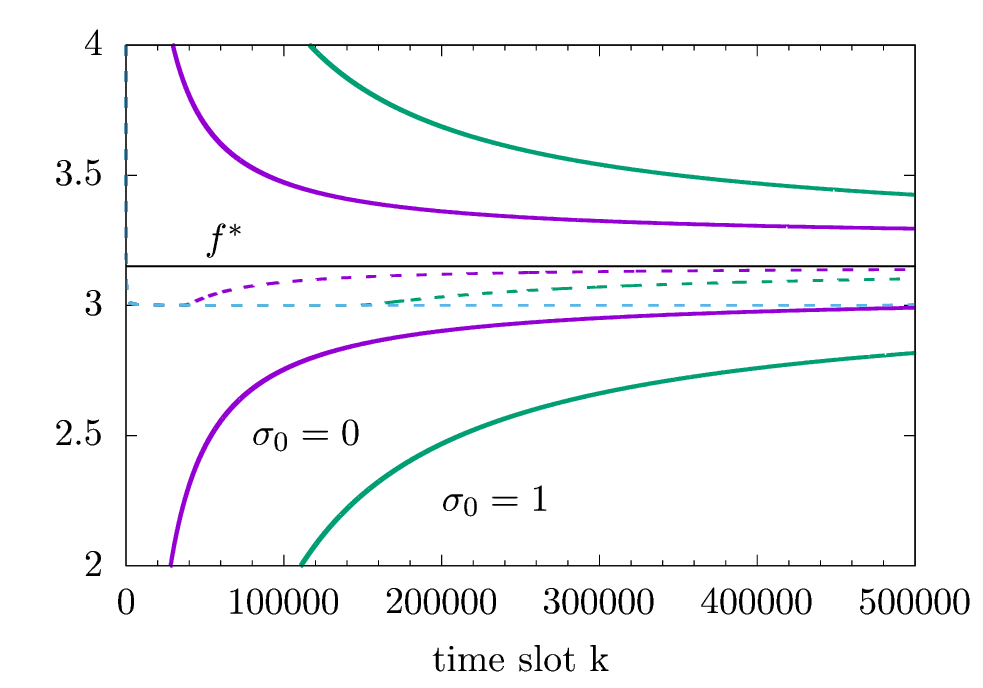}
\caption{Illustrating the convergence of $f(\vv{z}^\diamond_k)$ to a ball around $f^*$ (straight line) of Example \ref{ex:bt3ex} when $\epsilon = 0.05$ and $\sigma_0 \in \{ 0,1,4 \}$. Dashed lines indicate $f(\vv{z}^\diamond_k)$ with $\bar k=81$ while thick lines indicate upper and lower bounds of Theorem \ref{th:maintheorem}.}
\label{fig:bounds}
\end{figure}

\begin{figure}
\centering
\includegraphics[width=\columnwidth]{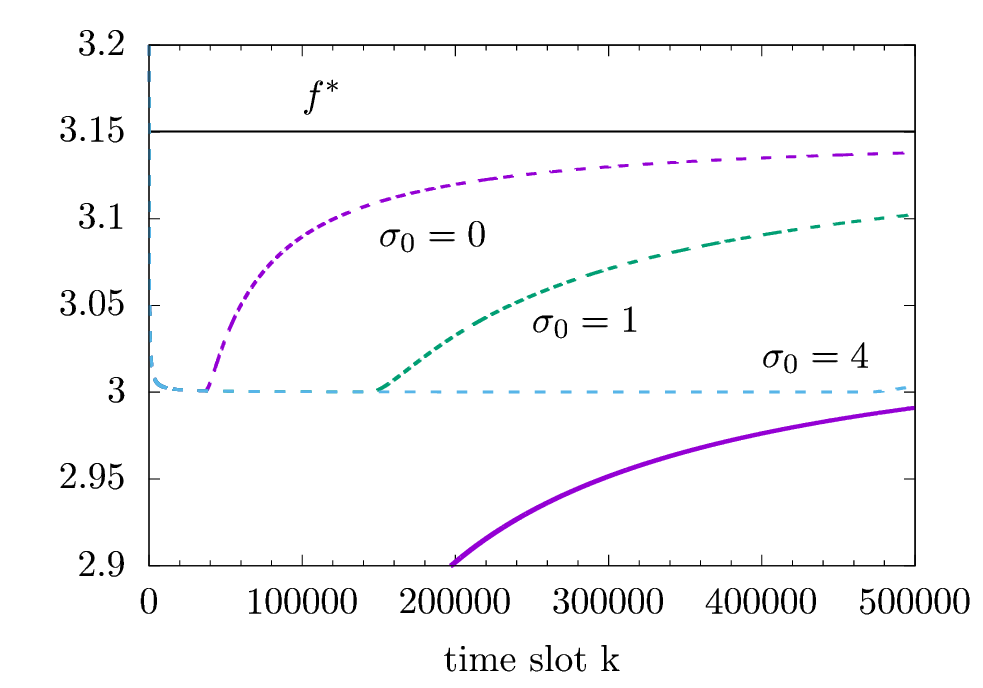}
\caption{Detail from Figure \ref{fig:bounds}.}
\label{fig:1zoom}
\end{figure}

Now suppose that $\vi{\tilde{\lambda}}{j} = \vi{\lambda}{j} + \alpha Y_k \sigma_0$ where $Y_k$ is uniformly randomly distributed between $-1$ and $1$. For $\sigma_0 \in \{0,1,4\}$ (so $\alpha \in \{7.29 \cdot 10^{-5} , 1.85 \cdot 10^{-5}  ,5.74 \cdot 10^{-6} \}$), Figure \ref{fig:bounds} plots $f(\vv{z}^\diamond_k)$ and the upper and lower bounds from Theorem \ref{th:maintheorem} vs $k$.  Figure \ref{fig:1zoom} shows detail from Figure \ref{fig:bounds}. It can be seen that, as expected, $f(\vv{z}^\diamond_k)$ is indeed upper and lower bounded by the values from Theorem \ref{th:maintheorem}.   It can also be seen that the upper and lower bounds are not tight, but they are not excessively loose either. 

\subsubsection{Violation of upper bound} 
\begin{figure}
\centering
\includegraphics[width=\columnwidth]{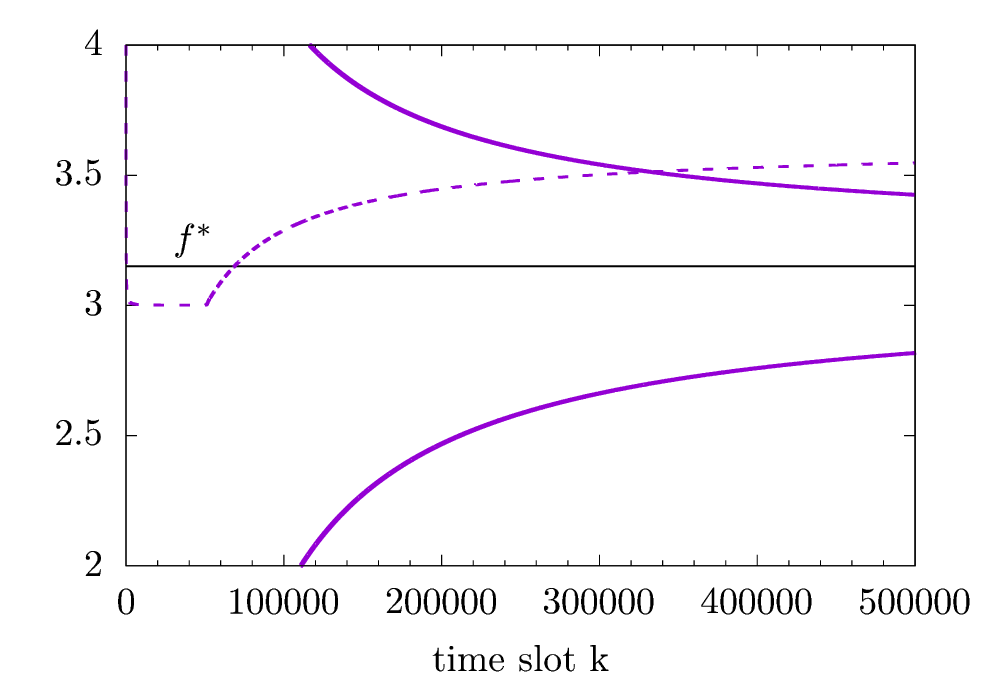}
\caption{Illustrating the violation of the bounds of Theorem \ref{th:maintheorem} when $\vi{\tilde{\lambda}}{j}_k = [\vi{{\lambda}}{j}_k +  \alpha e^k e^{-10^5}]^{[0, \bar \lambda]}$.  Dashed line indicates $f(\vv{z}^\diamond_k)$, $\bar k=  81$, while thicker lines indicate upper and lower bounds around $f^*$ (straight line).}
\label{fig:3v}
\end{figure}

Let $\vi{\tilde{\lambda}}{j}_k = [\vi{{\lambda}}{j}_k +  \alpha e^k e^{-10^5}]^{[0, \bar \lambda]}$.  With this choice the difference between $\vi{\lambda}{j}_k$ and $\vi{\tilde{\lambda}}{j}_k$ is uniformly bounded by $\alpha\sigma_0$ with $\sigma_0=1$  for $k \le 10^{5}$ but after that increases exponentially with $k$.   Figure \ref{fig:3v} plots $f(\vv{z}^\diamond_k)$ and the upper and lower bounds from Theorem \ref{th:maintheorem} when parameter $\alpha$ is selected according to Theorem \ref{th:maintheorem} assuming $\sigma_0 = 1$.   It can be seen that the upper and lower bounds hold  for $k \le 10^5$, but as the difference between multipliers increases $f(\vv{z}^\diamond_k)$ is not attracted to $f^*$ and it ends up violating the bounds.

\subsection{Example: Privacy-Enhancing Online Scheduling}
\label{ex:privacy}

We now present a small but interesting application example which illustrates some of the generality of Theorem \ref{th:maintheorem}.

Consider a sequence of information packets indexed by $k=1,2,\dots$.   Time is slotted and the packets arrive at a queue with inter-arrival times $\{b_k\}$, $k=1,2,\dots$ \emph{i.e.}, $b_k\in\mathbb{N}$ is the number of slots between the arrival of packet $k-1$ and packet $k$, with $b_1:=0$.   Outgoing packet $j$ is forwarded with inter-service time $s_j\in \{0,1,\dots,T\}\subset\mathbb{N}$ \emph{i.e.}, with $s_j$ slots between packet $x_j$ and the previously transmitted packet.    Dummy packets are transmitted as needed when no information packets are available, so as allow $s_j$ to be freely selected and to prevent large inter-arrivals times from propagating to the outgoing packet stream.  The aim is to select the queue service  $x_j$ such that the entropy of the inter-packet times of the outgoing packet stream is at least $E$, in order to provide a degree of resistence to traffic timing analysis, while stabilising the queue. 

The packet arrival process is not known in advance, other than the facts that it can be feasibily served, the inter-arrival times have finite mean $\lim_{k\rightarrow\infty}\frac{1}{k}\sum_{i=1}^k b_i = b$ and $|(\frac{1}{k}\sum_{i=1}^k b_i)-b|\le {\sigma_2}/k$ for some finite $\sigma_2>0$.

\begin{figure}
\centering
\includegraphics[width=\columnwidth]{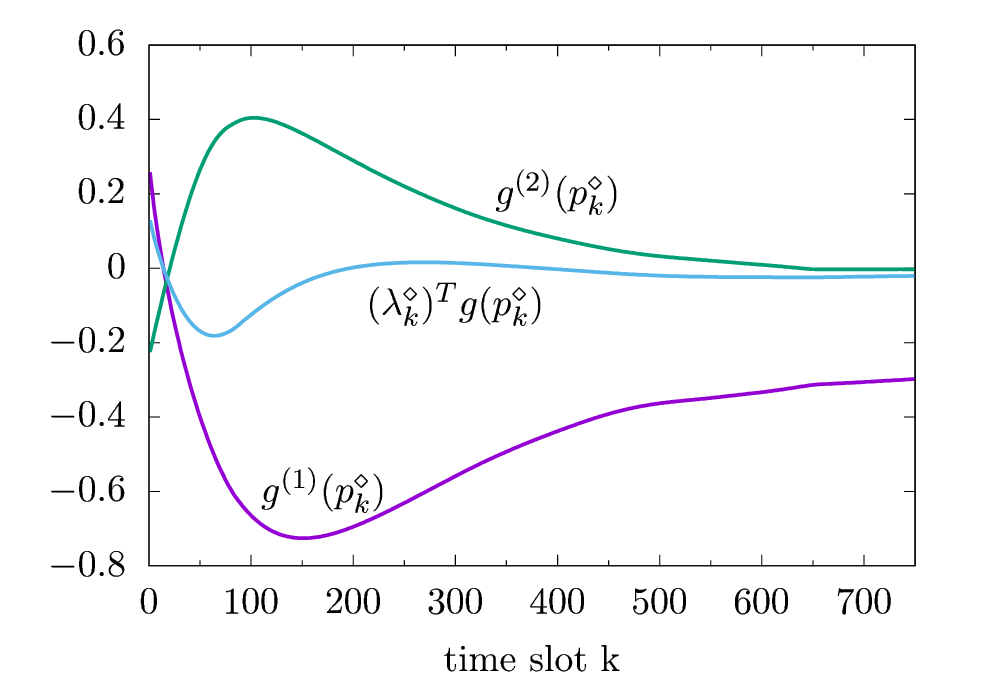}
\caption{Illustrating the convergence of sequence $\{ \vv{p}_k \}$ into the feasible region.  Updates (\ref{eq:priv_x_update}) - (\ref{eq:thetaupdate}) use parameters $T_\text{max} =5$, $E=\log(5)/5$, sequence $\{ b_k \} = \{0,1,0,1,0,\dots \}$ so $b = 1/2$, $\xi = b/2$, $\vi{\lambda}{1}_1 = \bar \lambda$, $\vi{\lambda}{2}_1 = 0$, $\bar\lambda = 1/2$ and $\alpha = \beta = 0.01$.}
\label{fig:privacy}
\end{figure}

Suppose the inter-service times $s_j$ are i.i.d. and let vector $\vv{p}$ with elements $\vi{p}{i}=\Prob(s_j=i)$, $i=0,\dots,T$ describe the probability distribution over set $\{0,1,\dots,T\}$.  The task can be formulated as the following feasibility problem (couched in convex optimisation form),
\begin{align*}
\min_{\vv{p}\in C} 1 \text{ s.t. } &\sum_{i=0}^{T_\text{max}} \vi{p}{i} \log \vi{p}{i} \le -E,\  \sum_{i=0}^{T_\text{max}} i \vi{p}{i} +\xi \le b
\end{align*}
where $\xi>0$ ensures that the mean inter-service time is strictly less than the mean inter-arrival time, so ensuring queue stability, and {$C:=\{\vv{p}\in [0,1]^T :\sum_{i=1}^T\vi{p}{i}\le 1\}$. } If the arrival process $\{b_k\}$ were known in advance, we could solve this optimisation to determine a feasible $\vv{p}$. When the arrivals are not known in advance, using generalised update (\ref{eq:gen1}) by Corollary \ref{th:corogenupdate} we can instead use the following online update to determine a sequence $\{ \vv{p}_k \}$ that converges to a feasible point. 
\begin{align}
 \vv{x}_k & \in \arg\min_{\vv{p}\in C} \vi{\lambda}{1}_k \vi{g}{1}(\vv{p}) + \vi{\lambda}{2}_k \vi{g}{2} (\vv{p}) \label{eq:priv_x_update}  \\
\vv{p}_{k+1}& =(1-\beta)\vv{p}_k+\beta\vv{x}_k  \\
\vi{\lambda}{1}_{k+1}& =\left[\vi{\lambda}{1}_k+\alpha \vi{g}{1} (\vv{p}_{k+1}) \right]^{[0,\bar \lambda]} \label{eq:phiupdate} \\
\vi{\lambda}{2}_{k+1}& =\left[\vi{\lambda}{2}_k+\alpha (\vi{g}{2} (\vv{p}_{k+1})+b-b_k ) \right]^{[0, \bar \lambda]}\label{eq:thetaupdate}
\end{align}
{{where $\vi{g}{1}(\vv{p}):= \sum_{i=0}^{T_\text{max}} \vi{p}{i} \log \vi{p}{i}+E $ and  $\vi{g}{2} (\vv{p}) := \sum_{i=0}^{T_\text{max}} i \vi{p}{i} + \xi -b$ with $\vi{\lambda}{1}_1,\vi{\lambda}{2}_1 \in [0,\bar \lambda]$.  The online update does not require knowledge of the mean  inter-arrival time $b$ since in (\ref{eq:priv_x_update}) the $\arg\min$ does not depend on $b$ while in (\ref{eq:thetaupdate}) we have  $\vi{g}{2}(\vv{p}_{k+1}) + b - b_k=\sum_{i=0}^{T_\text{max}} i \vi{p}{i}_{k+1} + \xi -b_k $.   

{\color{black}{Figure \ref{fig:privacy} illustrates the online update.  It can be seen that approximate complementary slackness converges to a ball around 0, and that constraints $\vi{g}{j} (\vv{p}^\diamond_k), \ j=1,2$ are attracted to the feasible region as $k$ increases.}}

}}
 
We highlight the following aspects of this example:

1) The online update differs from the standard dual-subgradient ascent in its use of the observed inter-arrival times $b_k$ rather than the (unknown) mean inter-arrival time $b$.  The inter-arrival times $b_k$ are discrete-valued, which also takes us outside of the usual setting.   The great advantage of the online update is that does not require knowledge of the mean rate $b$ of the packet arrival process, which is unknown beforehand, but only makes myopic use of available measurements to construct a packet schedule.

2) The constraint $\sum_{i=0}^{T_\text{max}} i \vi{p}{i} < b$ is expressed in terms of the packet inter-arrival and inter-service times rather than the number of packet arrivals and departures.   Hence, $\theta_k$ is not the scaled link queue occupancy but rather is related to the scaled queue waiting time.  Note that $\theta_k$ is not exactly the waiting time since the mean value $\sum_{i=0}^{T_\text{max}} i \vi{p}{i}$ is used for the inter-service time rather than the actual inter-service time realisations.

3) The transmission of dummy packets is explicitly an aspect of the application and it is the \emph{transmitted} packets (both dummy and information packets) which matter for satisfying the entropy constraint, not just the information packets.  This is because it is packet timing rather than packet content which is of interest here. 

4) Decision variable $x_k$ is not a packet or commodity scheduling action and so there are no issues around not having a packet to send when a queue is empty.

5) The entropy constraint is highly nonlinear and not at all like the type of flow constraint encountered in typical queueing network applications.

\section{Conclusions}

In this paper we investigate the connections between max-weight approaches and dual subgradient methods for convex optimisation.  We find that strong connections do indeed exist and we establish a clean, unifying theoretical framework that includes both max-weight and dual subgradient approaches as special cases.

\section{Appendix: Proofs}

\subsection{Proof of Theorems \ref{th:descent} and \ref{th:descent1b}}
The following two fundamental results are the key to establishing Theorem \ref{th:descent}:
\begin{lemma}\label{lem:two}
Let $D:=\{\vv{x}_1,\dots,\vv{x}_{|D|}\}$ be a finite set of points from $\mathbb R^n$ and $C := \conv{(D)}$. Then, for any point $\vv{y} \in C$ and vector $\vv{z} \in \mathbb R^n$ there exists a point $\vv{x} \in D$ such that $\vv{z}^T(\vv{x} - \vv{y}) \le 0$. 
\end{lemma}
\begin{IEEEproof}
Since $\vv{y}\in C:=\conv(D)$, $\vv{y}=\sum_{j=1}^{|D|}\vi{\theta}{j} \vv{x}_j$ with $\sum_{j=1}^{|D|} \vi{\theta}{j}=1$ and $\vi{\theta}{j}\in[0,1]$.  Hence, $\vv{z}^T(\vv{x}-\vv{y})=\sum_{j=1}^{|D|} \vi{\theta}{j}\vv{z}^T(\vv{x}-\vv{x}_j)$.  Select $\vv{x} \in \arg \min_{\vv{w}\in D} \vv{z}^T\vv{w}$.     Then $\vv{z}^T\vv{x} \le \vv{z}^T\vv{x}_j$ for all $\vv{x}_j \in D$ and so $\vv{z}^T(\vv{x}-\vv{y})\le 0$.
\end{IEEEproof}

\begin{lemma}[Non-Convex Descent]\label{th:decrease}
Let $F(\vv{z})$ be a convex function and suppose points $\vv{y}$, $\vv{z}\in C=\conv(D)$ exist such that $F(\vv{y})\le F(\vv{z})-\epsilon$, $\epsilon>0$.  Suppose $F(\cdot)$ has bounded curvature on $C$ with curvature constant $\mu_F$.  Then there exists at least one $\vv{x}\in D$ such that $F((1-\beta) \vv{z} + \beta \vv{x}) \le F(\vv{z})-\gamma\beta\epsilon$ with $\gamma \in (0,1)$ provided $0 < \beta\le(1-\gamma)\min\{\epsilon / (\mu_F\bar{x}^2_D),1\}$.
\end{lemma}
\begin{IEEEproof}
By convexity,
\begin{align*}
F(\vv{z}) +  \partial F(\vv{z})^T(\vv{y}-\vv{z}) \le F(\vv{y}) \le F(\vv{z})-\epsilon.
\end{align*}
Hence, $\partial F(\vv{z})^T(\vv{y}-\vv{z}) \le -\epsilon$.  Now observe that for $\vv{x}\in D$ we have $(1-\beta) \vv{z} + \beta \vv{x}\in C$ and by the bounded curvature of $F(\cdot)$ on $C$ 
\small
\begin{align*}
&F((1-\beta) \vv{z} + \beta \vv{x}) \\
&\le F(\vv{z}) + \beta\partial F(\vv{z})^T(\vv{x}-\vv{z}) + \mu_F \beta^2\|\vv{x}-\vv{z}\|_2^2\\
&= F(\vv{z}) + \beta\partial F(\vv{z})^T(\vv{y}-\vv{z})+\beta\partial F(\vv{z})^T(\vv{x}-\vv{y}) + \mu_F\beta^2\|\vv{x}-\vv{z}\|_2^2\\
&\le F(\vv{z}) - \beta\epsilon +\beta\partial F(\vv{z})^T(\vv{x}-\vv{y}) + \mu_F\beta^2\|\vv{x}-\vv{z}\|_2^2
\end{align*}
\normalsize
By Lemma \ref{lem:two} we can select $\vv{x}\in D$ such that $\partial F(\vv{z})^T(\vv{x}-\vv{y})\le 0$.  With this choice of $\vv{x}$ it follows that
\begin{align}
F((1-\beta) \vv{z} + \beta \vv{x}) 
&\le F(\vv{z}) - \beta\epsilon  + \mu_F\beta^2\|\vv{x}-\vv{z}\|_2^2 \notag \\
& \le F(\vv{z}) - \beta\epsilon  + \mu_F\beta^2\bar{x}_D^2 \label{eq:thdistanceup}
\end{align}
where (\ref{eq:thdistanceup}) follows from Lemma \ref{lem:one}, and the result now follows.
\end{IEEEproof}

\begin{IEEEproof}[Proof of Theorem \ref{th:descent}]{Since $F_k(\cdot)$ has bounded curvature for any $k$ it is continuous, and as $C$ is closed and bounded we have by the {Weierstrass theorem} (\emph{e.g.,} see Proposition 2.1.1 in \cite{convexanalysis}) that  $\min_{\vv{z} \in C} F_{k}(\vv{z})$ is finite.} We now proceed considering two cases: 

Case (i): $F_k(\vv{z}_k) - F_k(\vv{y}^*_k) \ge \epsilon$. By Lemma \ref{th:decrease} there exists $\vv{x}_k \in D$ such that $F_k ((1-\beta)\vv{z}_{k} + \beta \vv{x}_k) - F_k(\vv{z}_k) = F_k(\vv{z}_{k+1}) - F_k(\vv{z}_k) \le - \gamma\beta \epsilon$. Further, since $F_{k+1}(\vv{z}_{k+1}) - F_{k}(\vv{z}_{k+1}) \le \gamma_1 \gamma\beta \epsilon$ and $F_k(\vv{z}_k)  - F_{k+1} (\vv{z}_k) \le  \gamma_1 \gamma \beta \epsilon$ it follows
\begin{align}
F_{k+1}(\vv{z}_{k+1}) - F_{k+1}(\vv{z}_k) &  \le 2\gamma_1 \gamma \beta\epsilon - \gamma \beta\epsilon < 0.
\end{align}
That is, $F_k(\cdot)$ and $F_{k+1}(\cdot)$ decrease monotonically when $F_k(\vv{z}_k)- F_k(\vv{y}^*_k) \ge \epsilon$.

Case (ii): $F_k(\vv{z}_k) - F_k(\vv{y}^*_k) < \epsilon$. It follows that $F_k(\vv{z}_k) < F_k(\vv{y}_k^*)+\epsilon$.  Since $F_k(\cdot)$ is convex and has bounded curvature, $F_k(\vv{z}_{k+1})\le F_k(\vv{z}_k) + \beta\partial F_k(\vv{z}_k)^T(\vv{x}_{k}-\vv{z}_k)+\beta^2\mu_F\bar{x}_D^2 \le  F_k(\vv{y}_k^*)+\epsilon + \beta\partial F_k(\vv{z}_k)^T(\vv{x}_{k}-\vv{z}_k)+\beta^2\mu_F\bar{x}_D^2$.  The final term holds uniformly for all $\vv{x}_k\in D$ and since we select $\vv{x}_k$ to minimise $F_k(\vv{z}_{k+1})$ by Lemma \ref{lem:two} we therefore have $F_k(\vv{z}_{k+1}) \le F_k(\vv{y}^*_k)+\epsilon +\beta^2\mu_F\bar{x}_D^2$. Using the stated choice of $\beta$ and the fact that $F_{k+1}(\vv{z}_{k+1}) - \gamma_1 \gamma \beta \epsilon \le F_k(\vv{z}_{k+1}) $ yields
\begin{align}
F_{k+1}(\vv{z}_{k+1}) - F_{k}(\vv{y}^*_k)   \le \gamma_1 \gamma \beta\epsilon  + \epsilon + \beta (1-\gamma)\epsilon.
\end{align}Finally, since $F_k(\vv{y}^*_k) \le  F_k(\vv{y}_{k+1}^*) \le F_{k+1}(\vv{y}_{k+1}^*) +\gamma_1 \gamma \beta \epsilon$  we obtain
\begin{align*}
F_{k+1}(\vv{z}_{k+1})-F_{k+1}(\vv{y}_{k+1}^*) & \le 2\gamma_1 \gamma \beta\epsilon + \epsilon +  \beta (1-\gamma)\epsilon \\
& \le 2 \epsilon.
\end{align*}

We therefore have that $F_{k+1}(\vv{z}_k)$ is strictly decreasing when $F_k(\vv{z}_k)- F_k(\vv{y}^*_k) \ge \epsilon$ and otherwise uniformly upper bounded by $2\epsilon$.  It follows that for all $k$ sufficiently large $F_{k}(\vv{z}_{k+1})-F_{k}(\vv{y}_{k}^*) \le 2\epsilon$ as claimed.
\end{IEEEproof}

\begin{IEEEproof}[Proof of Theorem \ref{th:descent1b}]
Firstly, we make the following observations,
\begin{align}
&\arg\min_{\vv{z}\in C} F_k(\vv{z}_k) + \partial F_k(\vv{z}_k)^T(\vv{z}-\vv{z}_k)\notag\\
&\stackrel{(a)}{=}\arg\min_{\vv{z}\in C} \partial F_k(\vv{z}_k)^T\vv{z}
\stackrel{(b)}=\arg\min_{\vv{x}\in D} \partial F_k(\vv{z}_k)^T\vv{x}
\end{align}   
where equality $(a)$ follows by dropping terms not involving $\vv{z}$ and $(b)$ from the observation that we have a linear program (the objective is linear and set $C$ is a polytope, so defined by linear constraints) and so the optimum lies at an extreme point of set $C$ \emph{i.e.,} in set $D$.  We also have that
\begin{align*}
\small F_k(\vv{z}_k) + \partial F_k(\vv{z}_k)^T(\vv{x}_{k}-\vv{z}_k)
& \small \stackrel{(a)}{\le} F_k(\vv{z}_k) + \partial F_k(\vv{z}_k)^T(\vv{y}^*_k-\vv{z}_k)\\
& \small \stackrel{(b)}{\le} F_k(\vv{y}^*_k)
\le F_k(\vv{z}_k)
\end{align*}
\normalsize
where $\vv{y}^*_k\in\arg\min_{\vv{z}\in C}F_k(\vv{z})$, inequality $(a)$ follows from the minimality of $\vv{x}_k$ in $C$ noted above and $(b)$ from the convexity of $F_k(\cdot)$.   It follows that $\partial F_k(\vv{z}_k)^T(\vv{x}_{k}-\vv{z}_k)\le - (F_k(\vv{z}_k)  - F_k(\vv{y}^*_k)) \le 0$.   We have two cases to consider.  
Case (i):  $F_k(\vv{z}_k)-F_k(\vv{y}^*_k)\ge \epsilon$.  By the bounded curvature of $F_k(\cdot)$,
\begin{align*}
F_k(\vv{z}_{k+1}) 
&\le F_k(\vv{z}_k)  + \beta\partial F_k(\vv{z}_k)^T(\vv{x}_{k}-\vv{z}_k)+\mu_F\beta^2\bar{x}_D\\
&\le  F_k(\vv{z}_k) -\beta\epsilon+\mu_f\beta^2\bar{x}_D
\le  F_k(\vv{z}_k) -\gamma\beta\epsilon.
\end{align*}
Hence, 
\begin{align*}
F_{k+1}(\vv{z}_{k+1}) 
&\le F_k(\vv{z}_{k+1})  + |F_{k+1}(\vv{z}_{k+1})-F_k(\vv{z}_{k+1})|\\
 &\le F_k(\vv{z}_k) -\gamma\beta\epsilon  +\gamma_1\gamma\beta\epsilon,
\end{align*}and since $F_k(\vv{z}_k) \le F_{k+1}(\vv{z}_k) + \gamma_1 \gamma \beta  \epsilon$ we have that $F_{k+1}(\vv{z}_{k+1}) - F_k(\vv{z}_{k+1}) < 0$.
Case (ii): $F_k(\vv{z}_k)-F_k(\vv{y}^*_k)< \epsilon$.  Then 
\begin{align*}
F_k(\vv{z}_{k+1}) 
&\le F_k(\vv{z}_k)  + \beta\partial F_k(\vv{z}_k)^T(\vv{x}_{k}-\vv{z}_k)+\mu_F\beta^2\bar{x}_D\\
&\le F_k(\vv{y}^*_k ) + \epsilon + \beta \epsilon,
\end{align*}
and similar to the proof of Theorem \ref{th:descent} we obtain that $F_{k+1}(\vv{z}_{k+1}) - F_{k+1}(\vv{y}^*_{k+1}) \le 2\epsilon$.
We therefore have that $F_k(\vv{z}_{k})$ is strictly decreasing when  $F_k(\vv{z}_k)-F_k(\vv{y}^*_k)\ge \epsilon$ and otherwise uniformly upper bounded by $2\epsilon$.  Thus for $k$ sufficiently large $F_k(\vv{z}_{k+1})-F_k(\vv{y}^*_k)\le 2\epsilon$.
\end{IEEEproof}

\begin{IEEEproof}[Proof of Lemma \ref{th:subgradient}]
{\color{black}{
Let $\vv{\theta} \in \mathbb R^m$ such that $\vi{\theta}{j} \le \bar \lambda$ for all $j=1,\dots,m$ and see that
\begin{align}
& \| \vv{\lambda}_{k+1}-\vv{\theta}\|_2^2  \notag \\
& = \|[\vv{\lambda}_k + \alpha \vv{g}(\vv{z}_{k+1})]^{[0,\bar{\lambda}]}-\vv{\theta}\|_2^2 \notag \\
& \le \|[\vv{\lambda}_k + \alpha \vv{g}(\vv{z}_{k+1})]^{+}-\vv{\theta}\|_2^2 \label{eq:low1} \\
&  \le \| \vv{\lambda}_k + \alpha \vv{g}(\vv{z}_{k+1})-\vv{\theta}\|_2^2  \notag \\
& = \| \vv{\lambda}_k-\vv{\theta}\|_2^2 + 2\alpha(\vv{\lambda}_k-\vv{\theta})^T\vv{g}(\vv{z}_{k+1}) + \alpha^2 \| \vv{g}(\vv{z}_{k+1}) \|_2^2 \notag \\
& \le \| \vv{\lambda}_k-\vv{\theta}\|_2^2 + 2\alpha(\vv{\lambda}_k-\vv{\theta})^T\vv{g}(\vv{z}_{k+1}) + \alpha^2 m {\bar g}^2, \label{eq:low2}
\end{align} where $(\ref{eq:low1})$ follows since $\bar \lambda \ge \vi{\theta}{j}$ and $(\ref{eq:low2})$ from the fact that $\|\vv{g}(\vv{z})\|^2_2 \le m{\bar g}^2$ for all $\vv{z} \in C$. Applying the latter argument recursively for $i=1, \dots, k$ yields $\|\vv{\lambda}_{k+1}-\vv{\theta}\|_2^2 \le \| \vv{\lambda}_{1}-\vv{\theta}\|_2^2 + 2\alpha \sum_{i = 1}^{k} (\vv{\lambda}_i-\vv{\theta})^T\vv{g}(\vv{z}_{i+1}) + \alpha^2 m {\bar g}^2 k $. Rearranging terms, dividing by $2 \alpha k$, and using the fact that  $\|\vv{\lambda}_{k+1} - \vv{\theta}\|_2^2 \ge 0$ and $\| \vv{\lambda}_1 - \vv{\theta} \|_2^2 \le 2 m \bar \lambda^2$ we have
\begin{align}
- \frac{m \bar \lambda^2}{\alpha k } -  \frac{\alpha}{2}m {\bar g}^2
& \le \frac{1}{k} \sum_{i=1}^{k} (\vv{\lambda}_i-\vv{\theta})^T\vv{g}(\vv{z}_{i+1}) \\
& = \frac{1}{k} \sum_{i=1}^{k} L(\vv{z}_{i+1},\vv{\lambda}_i)  - L(\vv{z}_{i+1},\vv{\theta}).
 \label{eq:subgradup}
\end{align}
Next, see that by the definition of sequence $\{ \vv{z}_k \}$ we can write $ \frac{1}{k} \sum_{i=1}^k L(\vv{z}_{i+1}, \vv{\lambda}_i) \le \frac{1}{k} \sum_{i=1}^k  q(\vv{\lambda}_i) + 2 \epsilon \le q(\vv{\lambda}^\diamond_k) + 2 \epsilon$ where the last inequality follows by the concavity of $q$. That is, 
\begin{align}
- \frac{m \bar \lambda^2}{\alpha k } -  \frac{\alpha}{2}m{\bar g}^2 - 2\epsilon
& \le q(\vv{\lambda}^\diamond_k) - \frac{1}{k} \sum_{i=1}^k L(\vv{z}_{i+1},\vv{\theta})
\end{align}By fixing $\vv{\theta}$ to $\vv{\lambda}^*$ and $\vv{\lambda}^\diamond_k$ and using the fact that $ \frac{1}{k} \sum_{i=1}^k L(\vv{z}_{i+1},\vv{\lambda}^\diamond_k) \ge  L(\vv{z}^\diamond_k, \vv{\lambda}^\diamond_k)$ for all $k=1,2,\dots$ and $ \frac{1}{k} \sum_{i=1}^k L(\vv{z}_{i+1},\vv{\lambda}^*) \ge f^*$ we have that
\begin{align}
- \frac{m \bar \lambda^2}{\alpha k } -  \frac{\alpha}{2}m {\bar g}^2
- 2 \epsilon & \le q(\vv{\lambda}^\diamond_k) - f^* \le 0 \label{eq:avgup}
\end{align}and
\begin{align}
- \frac{m \bar \lambda^2}{\alpha k } -  \frac{\alpha}{2}m {\bar g}^2
- 2 \epsilon & \le q(\vv{\lambda}^\diamond_k) - L(\vv{z}^\diamond_k, \vv{\lambda}^\diamond_k) \le 0  \label{eq:avglow}.
\end{align}Multiplying (\ref{eq:avgup}) by $-1$ and combining it with (\ref{eq:avglow}) yields the result.
}}
\end{IEEEproof}

\begin{IEEEproof}[Proof of Lemma \ref{th:slackness}]

We start by showing that updates $[\vv{\lambda}_k + \alpha \vv{g}(\vv{z}_{k+1})]^+$ and $[\vv{\lambda}_k + \alpha \vv{g}(\vv{z}_{k+1})]^{[0, \bar \lambda]}$ are interchangeable when $L(\vv{z}_{k+1}, \tilde{\vv{\lambda}}_k)$ is uniformly close to $q(\tilde{\vv{\lambda}}_k)$. First of all see that 
\begin{align*}
& \| \vv{\lambda}_{k+1} - \vv{\lambda}^* \|_2^2  \\
& = \| [\vv{\lambda}_k + \alpha \vv{g}(\vv{z}_{k+1})]^+ - \vv{\lambda}^* \|_2^2 \\
& \le \| \vv{\lambda}_k + \alpha \vv{g}(\vv{z}_{k+1}) - \vv{\lambda}^* \|_2^2 \\
& = \| \vv{\lambda}_k - \vv{\lambda}^* \|_2^2 + \alpha^2 \| \vv{g}(\vv{z}_{k+1}) \|_2^2  + 2 \alpha (\vv{\lambda}_k - \vv{\lambda}^*)^T \vv{g}(\vv{z}_{k+1}) \\
& \le \| \vv{\lambda}_k - \vv{\lambda}^* \|_2^2 + \alpha^2 m {\bar g}^2  + 2 \alpha (\vv{\lambda}_k - \vv{\lambda}^*)^T \vv{g}(\vv{z}_{k+1}).
\end{align*}Now observe that since $\| \vv{\lambda}_k - \tilde{\vv{\lambda}}_k \|_2 \le \| \vv{\lambda}_k - \tilde{\vv{\lambda}}_k \|_1 \le \alpha m \sigma_0 $ we can write $  (\vv{\lambda}_k - \vv{\lambda}^*)^T \vv{g}(\vv{z}_{k+1}) =   (\tilde{\vv{\lambda}}_k - \vv{\lambda}^*)^T \vv{g}(\vv{z}_{k+1}) +  (\vv{\lambda}_k - \tilde{\vv{\lambda}}_k)^T \vv{g}(\vv{z}_{k+1}) \le   (\tilde{\vv{\lambda}}_k - \vv{\lambda}^*)^T \vv{g}(\vv{z}_{k+1}) +  \| \vv{\lambda}_k - \tilde{\vv{\lambda}}_k\|_2  \| \vv{g}(\vv{z}_{k+1}) \|_2 \le  (\tilde{\vv{\lambda}}_k - \vv{\lambda}^*)^T \vv{g}(\vv{z}_{k+1}) +  \alpha m^2 \sigma_0 {\bar g} =  L(\vv{z}_{k+1}, \tilde{\vv{\lambda}}_k) - L(\vv{z}_{k+1},\vv{\lambda}^*) +  \alpha m^2 \sigma_0 {\bar g}$. Furthermore, since $L(\vv{z}_{k+1}, \tilde{\vv{\lambda}}_k) \le q(\tilde{\vv{\lambda}}_k) + 2\epsilon$ and $- L(\vv{z}_{k+1},\vv{\lambda}^*) \le -q(\vv{\lambda}^*)$ it follows that
\begin{align}
& \| \vv{\lambda}_{k+1} - \vv{\lambda}^* \|_2^2  - \| \vv{\lambda}_k - \vv{\lambda}^* \|_2^2 \notag \\
&  \quad  \le  \alpha^2  (m{\bar g}^2 + 2m^2  \sigma_0 {\bar g})  + 2 \alpha ( q(\tilde{\vv{\lambda}}_k) + 2\epsilon - q(\vv{\lambda}^*)). \label{eq:distmulti}
\end{align}

Now let $Q_\delta : = \{ {\vv{\lambda}} \succeq \vv{0} : q({\vv{\lambda}}) \ge q(\vv{\lambda}^*) - \delta \}$ and consider two cases. Case (i) $(\tilde{\vv{\lambda}}_k \notin Q_{\delta})$. Then $q(\tilde{\vv{\lambda}}_k) - q(\vv{\lambda}^*) < - \delta$ and from (\ref{eq:distmulti}) we have that $\| \vv{\lambda}_{k+1} - \vv{\lambda}^* \|_2^2 < \| \vv{\lambda}_{k} - \vv{\lambda}^* \|_2^2$, \emph{i.e.,} 
 \begin{align*}
 \| \vv{\lambda}_{k+1} - \vv{\lambda}^* \|_2  - \| \vv{\lambda}_k - \vv{\lambda}^* \|_2 < 0
 \end{align*} and so $\vv{\lambda}_k$ converges to a ball around $\vv{\lambda}^*$ when $\tilde{\vv{\lambda}}_k \in Q_\delta$.  Case (ii) $(\tilde{\vv{\lambda}}_k \in Q_{\delta})$. See that  $\| \vv{\lambda}_{k+1} - \vv{\lambda}^* \| = \| [\vv{\lambda}_k + \alpha \vv{g}(\vv{z}_{k+1})]^+ - \vv{\lambda}^* \|_2 \le \| \vv{\lambda}_k + \alpha \vv{g}(\vv{z}_{k+1}) - \vv{\lambda}^* \|_2 \le \| \vv{\lambda}_k \|_2 + \| \vv{\lambda}^* \|_2 + \alpha m {\bar g}$. Next recall that when the Slater condition holds by Lemma \ref{th:setq} we have for all $\vv{\lambda} \in Q_\delta$ then $\| \vv{\lambda} \| \le \frac{1}{\upsilon}(\eta + \delta)$ where $\eta:=f(\bar{\vv{z}}) - q(\vv{\lambda}^*)$ and $\bar{\vv{z}}$ a Slater vector. Therefore, 
 \begin{align*}
 \| \vv{\lambda}_{k+1} - \vv{\lambda}^* \|_2 \le \frac{2}{\upsilon} (\eta + \delta) + \alpha m {\bar g}.
 \end{align*}
 
From both cases it follows that if $ \| \vv{\lambda}_1 - \vv{\lambda}^* \|_2 \le \frac{2}{\upsilon} (\eta + \delta) + \alpha m {\bar g}$ then $\|\vv{\lambda}_k - \vv{\lambda}^* \|_2 \le \frac{2}{\upsilon} (\eta + \delta) + \alpha m {\bar g}$ for all $k \ge 1$. Using this observation and the fact that $\|\vv{\lambda}_1 - \vv{\lambda}^* \|_2 \ge | \| \vv{\lambda}_1\|_2 - \| \vv{\lambda}^*\|_2 | \ge   \| \vv{\lambda}_1\|_2 - \| \vv{\lambda}^*\|_2 $ we obtain that when $\| \vv{\lambda}_1 \|_2 \le \frac{3}{\upsilon}(\eta + \delta) + \alpha m {\bar g}$ then $\| \vv{\lambda}_k \|_2 \le  \frac{3}{\upsilon}(\eta + \delta) + \alpha m {\bar g}$ for all $k \ge 1$. That is, if we choose $ \vi{\lambda}{j}_1 \le \frac{3}{ \upsilon}(\eta + \delta)  + \alpha m {\bar g} \le \bar \lambda $ then $\vi{\lambda}{j}_k \le \bar \lambda$ for all $j=1,\dots,m$, $k \ge 1$ and so updates $[\vv{\lambda}_k + \vv{g}(\vv{z}_{k+1})]^+$ and $[\vv{\lambda}_k + \vv{g}(\vv{z}_{k+1})]^{[0,\bar \lambda]}$ are interchangeable as claimed. 

Now we proceed to prove the upper and lower bounds in (\ref{eq:slackness}). For the lower bound see first that
\begin{align*}
 \| \vv{\lambda}_{k+1} \|^2_2  
&  = \| [\vv{\lambda}_k + \alpha \vv{g}(\vv{z}_{k+1})]^+ \|^2_2 \\
&  \le \| \vv{\lambda}_k  + \alpha \vv{g}(\vv{z}_{k+1})\|^2_2 \\
&  = \| \vv{\lambda}_k \|^2_2 + \alpha^2 \| \vv{g}(\vv{z}_{k+1}) \|^2_2 + 2 \alpha \vv{\lambda}^T_k \vv{g}(\vv{z}_{k+1}) \\
&  \le \| \vv{\lambda}_k \|^2_2 + \alpha^2 m {\bar g}^2 + 2 \alpha \vv{\lambda}^T_k \vv{g}(\vv{z}_{k+1}) 
\end{align*}
Rearranging terms and applying the latter bound recursively for $i=1,\dots,k$ yields $ 2 \alpha \sum_{i=1}^k  {\vv{\lambda}}^T_i \vv{g}(\vv{z}_{i+1})  \ge  \| \vv{\lambda}_{k+1} \|^2_2 - \| \vv{\lambda}_1\|^2_2 - \alpha^2 m {\bar g}^2 k \ge - \| \vv{\lambda}_1\|^2_2 - \alpha^2 m {\bar g}^2 k $. {\color{black}{The bound does not depend on sequence $\{ \vv{z}_k \}$, hence, it holds for any sequence of points in $C$. Fixing $\vv{z}_{i+1}$ to $\vv{z}^\diamond_k$ for all $i=1,\dots,k$ we can write $2 \alpha \sum_{i=1}^k  {\vv{\lambda}}^T_i \vv{g}(\vv{z}^\diamond_{k}) = 2 \alpha k (\vv{\lambda}^\diamond_k)^T \vv{g}(\vv{z}^\diamond_{k})$. Dividing by $2\alpha k $ and using the fact that $\| \vv{\lambda}_1\|_2^2 \le m \bar \lambda^2$ yields
\begin{align*}
- \frac{m \bar \lambda}{2\alpha k} - \frac{\alpha}{2} m {\bar g}^2 \le  (\vv{\lambda}^\diamond_k)^T \vv{g}(\vv{z}^\diamond_{k}).
\end{align*}
For the upper bound see that
$\vv{\lambda}_{k+1} = [\vv{\lambda}_k + \alpha \vv{g}(\vv{z}_{k+1})]^+ \succeq \vv{\lambda}_{k} + \alpha \vv{g}(\vv{z}_{k+1})$ 
and so we can write
$\alpha \sum_{i=1}^{k} \vv{g}(\vv{z}_{i+1}) \preceq \sum_{i=1}^k (\vv{\lambda}_{i+1}-\vv{\lambda}_i) = \vv{\lambda}_{k+1} - \vv{\lambda}_1 \preceq \vv{\lambda}_{k+1} $. 
Next, by the convexity of $\vv{g}(\cdot)$ we have that 
$\frac{1}{k} \sum_{i=1}^{k} \alpha \vv{g}(\vv{z}_{i+1}) \succeq \alpha \vv{g}(\vv{z}^\diamond_k)$ and so it follows that $\vv{g}(\vv{z}^\diamond_k) \preceq \vv{\lambda}_{k+1}/(\alpha k)$. Multiplying the last equation by $\vv{\lambda}^\diamond_k$ and using the fact that $\vv{0} \preceq\vv{\lambda}_{k+1} \preceq \bar\lambda \vv{1}$ and $\vv{0} \preceq \vv{\lambda}_k^\diamond \preceq \bar \lambda \vv{1}$ yields the upper bound.

Finally, the constraint violation bound (\ref{eq:feasibility}) follows from the fact that $\vv{g}(\vv{z}^\diamond_k) \preceq {\bar \lambda}^2/({\alpha k}) \vv{1}$. 
}}\end{IEEEproof}

\begin{IEEEproof}[Proof of Lemma \ref{th:sequences}]
First of all see that 
$
 | \lambda_{k+1} - \tilde{\lambda}_{k+1} |  
=  | [\lambda_k + \delta_k]^{[0,\bar \lambda]} - [\tilde \lambda_k + \tilde \delta_k]^{[0,\bar \lambda]} |
 \le  | [\lambda_k + \delta_k]^{[0,\bar \lambda]} - [\tilde \lambda_k + \tilde \delta_k]^+ | 
 =  | [\tilde \lambda_k + \tilde \delta_k]^+ -  [\lambda_k + \delta_k]^{[0,\bar \lambda]} | 
 \le  | [\tilde \lambda_k + \tilde \delta_k]^+ -  [\lambda_k + \delta_k]^+ |.
$
We now proceed to bound the RHS of the last equation. 
Let $\Delta_k := - \min(\lambda_k + \delta_k,0)$, \emph{i.e.,} $\lambda_{k+1}  = \lambda_k + \delta_k + \Delta_k$  so that we can write $\lambda_{k+1} = \lambda_1 + \sum_{i=1}^k (\delta_i + \Delta_i)$.  Note that when $\lambda_{k+1}> 0$ then $\Delta_k = 0$, and that when $\lambda_{k+1} = 0$ then $   \sum_{i=1}^k \Delta_i =  - \lambda_1 - \sum_{i=1}^k \delta_i$. Next, note that since $\Delta_k$ is nonnegative for all $k$ by construction we have that $\sum_{i=1}^k \Delta_i$ is non-decreasing in $k$. Using the latter observation it follows that $\sum_{i=1}^k \Delta_i = [- \lambda_1 - \min_{1 \le j \le k } \sum_{i=1}^j \delta_i]^+$ and therefore
\begin{align*}
\lambda_{k+1} =   \sum_{i=1}^k \delta_i + \max \left\{  \Theta_k , \lambda_1 \right\}
\end{align*}where $\Theta_k : = - \min_{1 \le j \le k} \sum_{i=1}^j \delta_i$. Now see that
\begin{align*}
& | \lambda_{k+1} - \tilde \lambda_{k+1}| \\
& \textstyle =  |   \sum_{i=1}^k \delta_i + \max \{   \Theta_k , \lambda_1 \} -  \sum_{i=1}^k \tilde \delta_i - \max \{   \tilde \Theta_k , \lambda_1 \} |\\
& \textstyle \le  |   \sum_{i=1}^k \delta_i -  \tilde \delta_i | + | \max \{   \Theta_k , \lambda_1 \}  - \max \{   \tilde \Theta_k , \lambda_1 \} |\\
& \textstyle \stackrel{(a)}{\le}  |   \sum_{i=1}^k \delta_i - \tilde \delta_i | + |  \tilde \Theta_k - \Theta_k |\\
& \textstyle = | \sum_{i=1}^{k} \delta_i  - \tilde \delta_i | +    | \underset{1 \le j \le k}{\min} \sum_{i=1}^j  \tilde \delta_j  - \underset{1 \le j \le k}{\min} \sum_{i=1}^j  \delta_j  | \\
& \textstyle = | \sum_{i=1}^{k} \delta_i  - \tilde \delta_i | +    | \underset{1 \le j \le k}{\max} \sum_{i=1}^j  - \tilde \delta_j  - \underset{1 \le j \le k}{\max}\sum_{i=1}^j  -\delta_j  | \\
& \textstyle \le | \sum_{i=1}^{k} \delta_i  - \tilde \delta_i | +   \underset{1 \le j \le k}{\max} | \sum_{i=1}^j   \delta_j  - \sum_{i=1}^j \tilde \delta_j  |
\end{align*}where $(a)$ follows easily from enumerating the four cases. Finally, since $| \sum_{i=1}^{k} \delta_i  - \tilde \delta_i | \le \max_{i \le j \le k} | \sum_{i=1}^{j} \delta_i  - \tilde \delta_i |$ and $| \sum_{i=1}^{k} \delta_i  - \tilde \delta_i | \le \epsilon$  for all $k=1,2,\dots$ the result follows.
\end{IEEEproof}

\begin{IEEEproof}[Proof of Theorem \ref{th:auxiliarymultiplier}]
By Lemma \ref{th:sequences} we require $ |  \sum_{i=1}^k \vvi{a}{j}(\vv{z}_{i+1} - \vv{x}_i) +\vi{b}{j}_i-\vi{b}{j}|$ to be uniformly bounded in order to establish the boundedness of $ | \vi{\tilde{\lambda}}{j}_k- \vi{\lambda}{j}_k|$ for all $k\ge 1$. However, since $ |  \sum_{i=1}^k \vvi{a}{j}(\vv{z}_{i+1} - \vv{x}_i) +\vi{b}{j}_i-\vi{b}{j}| \le  |  \sum_{i=1}^k \vvi{a}{j}(\vv{z}_{i+1} - \vv{x}_i)| + | \sum_{i=1}^k\vi{b}{j}_i -\vi{b}{j}| $ and  $| \sum_{i=1}^k\vi{b}{j}_i -\vi{b}{j}| \le {\sigma_2}$ by assumption, it is sufficient to show that $|  \sum_{i=1}^k \vvi{a}{j}(\vv{z}_{i+1} - \vv{x}_i)|$ is bounded.

Now observe that since  $\vv{z}_{i+1} = (1-\beta) \vv{z}_{i} + \beta \vv{x}_{i}$ we have $\vv{z}_{i+1} - \vv{x}_i = (1-\beta)(\vv{z}_i - \vv{x}_i)$. That is, $\sum_{i=1}^k (\vv{z}_{i+1} - \vv{x}_i) = (1-\beta)\sum_{i=1}^k (\vv{z}_i - \vv{x}_i)$. Further, since $\sum_{i=1}^k (\vv{z}_i - \vv{x}_i) = \sum_{i=1}^{k-1} (\vv{z}_{i+1} - \vv{x}_i) + (\vv{z}_1 - \vv{x}_k) = (1-\beta)\sum_{i=1}^{k-1} (\vv{z}_i - \vv{x}_i) + (\vv{z}_1- \vv{x}_k)$ it follows that  $\sum_{i=1}^k (\vv{z}_{i+1} - \vv{x}_i) = (1-\beta)^2 \sum_{i=1}^{k-1}(\vv{z}_{i} - \vv{x}_i) + (1-\beta) (\vv{z}_1 - \vv{x}_k)$.  Applying the preceding argument recursively we obtain that $\sum_{i=1}^k (\vv{z}_{i+1} - \vv{x}_i) =  (1-\beta)(\vv{z}_1 - \vv{x}_k) + (1-\beta)^2(\vv{z}_1 - \vv{x}_{k-1}) + \dots + (1-\beta)^{k} (\vv{z}_1 - \vv{x}_{1})$, \emph{i.e.,} 
\begin{align}
\sum_{i=1}^k (\vv{z}_{i+1} - \vv{x}_i)  = \sum_{i=1}^{k} (1-\beta)^{k+1-i} (\vv{z}_1 - \vv{x}_i). \label{eq:sumexpandedfinal}
\end{align}
Using (\ref{eq:sumexpandedfinal}) it follows that 
\begin{align}
& \textstyle 2 \alpha \left|  \sum_{i=1}^k \vvi{a}{j}(\vv{z}_{i+1} - \vv{x}_i) \right| \notag \\
& \textstyle \quad \le 2 \alpha \left|  \sum_{i=1}^k (1-\beta)^{k+1-i} \vvi{a}{j}(\vv{z}_{1}-{\vv{x}}_i) \right|  \notag \\
& \textstyle \quad \le 2 \alpha{\sigma_1}   \sum_{i=1}^k (1-\beta)^{k+1-i}  \label{eq:zkupdate}
\end{align}
where ${\sigma_1} :=  2 \max_{\vv{z} \in C} \| \vv{A} \vv{z} \|_\infty$. Next, see that $\sum_{i=1}^k (1-\beta)^{k+1-i} = (1-\beta)^{k+1} \sum_{i=1}^k (1-\beta)^{-i}$ and that
\begin{align*}
\sum_{i=1}^k \frac{1}{(1-\beta)^i} = \frac{1-(1-\beta)^{k+1}}{\beta(1-\beta)^k} .
\end{align*}Therefore, 
$\sum_{i=1}^k (1-\beta)^{-i} <  (1-\beta)^{-k} / \beta$ and so 
\begin{align*}
(1-\beta)^{k+1}\sum_{i=1}^k (1-\beta)^{-i} < \frac{(1-\beta)}{\beta} < \frac{1}{\beta}.
\end{align*}Finally, using the latter bound in (\ref{eq:zkupdate}) the stated result now follows. 
\end{IEEEproof}

\bibliographystyle{IEEEtran}
\bibliography{IEEEabrv,bib}

\end{document}